\documentclass[11pt]{article}

\usepackage[T1]{fontenc}
\usepackage[latin1]{inputenc}
\usepackage{caption}
\usepackage{algorithm,algpseudocode}
\usepackage{graphicx,url}
\usepackage{amsmath,amssymb,amsopn}
\usepackage{tabularx}
\usepackage{color}
\usepackage{calc}
\usepackage{caption}
\usepackage{longtable}
\usepackage{multirow}
\usepackage{mathrsfs}
\usepackage{array}
\usepackage{hyperref}
\usepackage{tikz}
\usepackage{booktabs}
\usepackage{pgfplotstable}
\usepackage{pgfplots}
\usepackage[outline]{contour}
\usepackage{calrsfs}
\usepackage{graphicx}
\usepgfplotslibrary{groupplots}
\usetikzlibrary{matrix}
\usetikzlibrary{automata, external, shapes.geometric, fit}
\def \bv{\vec{b}}

\def \dv{\vec{d}}

\def \fv{\vec{f}}
\def \gv{\vec{g}}

\def \pv{\vec{p}}
\def \qv{\vec{q}}
\def \rv{\vec{r}}

\def \uv{\vec{u}}
\def \vv{\vec{v}}
\def \wv{\vec{w}}
\def \xv{\vec{x}}
\def \yv{\vec{y}}
\def \zv{\vec{z}}

\def \um{\mat{u}}

\def \Am{\mat{A}}

\def \Cm{\mat{A_c}}

\def \Im{\mat{I}}

\def \Km{\mat{K}}

\def \Mm{\mat{M}}

\def \Qm{\mat{C}}
\def \QQm{\mat{Q}}
\def \Pm{\mat{P}}
\def \Rm{\mat{R}}

\def \Id{\mat{I}}

\def \R{\mathbb{R}}						\def \E {\mathbb{E}}

\DeclareMathAlphabet{\pazocal}{OMS}{zplm}{m}{n}

\renewcommand{\vec}[1]{\boldsymbol{#1}}
\newcommand{\mat}[1]{\boldsymbol{{#1}}}

\algnewcommand\algorithmiconput{\textbf{Constants:}}
\algnewcommand\algorithmicinput{\textbf{Input:}}
\algnewcommand\algorithmicoutput{\textbf{Output:}}
\algnewcommand{\algorithmicgoto}{\textbf{go to}}

\algnewcommand\Constants{\item[\algorithmiconput]}
\algnewcommand\Input{\item[\algorithmicinput]}\algnewcommand\Output{\item[\algorithmicoutput]}\algnewcommand{\Goto}[1]{\algorithmicgoto~\ref{#1}}

\definecolor{myblack}{RGB}{53, 53, 53}
\definecolor{myblue}{RGB}{40, 75, 99}
\definecolor{myred}{RGB}{192, 50, 33}
\definecolor{myyellow}{RGB}{255, 166, 48}
\definecolor{mywhite}{RGB}{240, 237, 238}
\definecolor{mygreen}{RGB}{0, 102, 0}

\definecolor{green1}{RGB}{9, 82, 86}
\definecolor{green2}{RGB}{8, 127, 140}
\definecolor{green3}{RGB}{6, 167, 125}
\definecolor{green4}{RGB}{79, 109, 122}
\definecolor{green5}{RGB}{192, 214, 223}
\definecolor{violet}{RGB}{26,69,131}

\definecolor{checkgreen}{rgb}{0,0.6,0}
\definecolor{phase1}{rgb}{0.008,0.655,1.000}
\definecolor{phase2}{rgb}{0.016,0.75,0.700}
\definecolor{phase3}{rgb}{0.929,0.35,0.700}
\definecolor{icsyellow}{cmyk}{0.00,0.11,0.53,0.00}

\definecolor{blackmy}{RGB}{38, 70, 83}
\definecolor{bluemy}{RGB}{39, 125, 161}
\definecolor{greenmy}{RGB}{42, 167, 143}
\definecolor{yellowmy}{RGB}{233, 196, 106}
\definecolor{brownmy}{RGB}{244, 162, 97}
\definecolor{redmy}{RGB}{249, 65, 68}

\definecolor{darkbluemy}{RGB}{65, 59, 147}
\definecolor{lightbluemy}{RGB}{71, 139, 194}
\definecolor{greenmy}{RGB}{98, 173, 153}
\definecolor{darkorangemy}{RGB}{230, 142, 52}
\definecolor{lightorangemy}{RGB}{217, 172, 59}

\definecolor{blue1}{RGB}{1, 42, 74}
\definecolor{blue2}{RGB}{1, 73, 124}
\definecolor{blue3}{RGB}{42, 111, 151}
\definecolor{blue4}{RGB}{44, 125, 160}
\definecolor{blue5}{RGB}{70, 143, 175}
\definecolor{blue6}{RGB}{137, 194, 217}

\definecolor{red1}{RGB}{204, 68, 75}
\definecolor{red2}{RGB}{218, 85, 82}
\definecolor{red3}{RGB}{227, 150, 149}
\definecolor{red4}{RGB}{228, 190, 171}

\definecolor{brown1}{RGB}{92,178,112}
\definecolor{brown2}{RGB}{130,194,110}		
\definecolor{brown3}{RGB}{163, 193, 173}

\usepackage[normalem]{ulem}
\usepackage{soul,xcolor}
\usepackage[many]{tcolorbox}
\usepackage{mathtools}
\usepackage{algorithm}
\usepackage{enumitem}
\usepackage{listings}
\usepackage{relsize}
\usepackage{nameref}

\usepackage{longtable}

\newcommand{\algorithmiccommentMine}[1]{\bgroup\hfill$\triangleright$~{ \textcolor{gray}{#1}}\egroup}
\newcommand\COMMENTmine[1]{\algorithmiccommentMine{#1}}

\newcommand\oldtext[1]{}
\newcommand\cancel[1]{}

\newcommand\oldtextt[1]{}

\usetikzlibrary{spy}
\usepackage{pgfplots}
\usepgfplotslibrary{fillbetween}

\usepackage{array}
\newcommand{\PreserveBackslash}[1]{\let\temp=\\#1\let\\=\temp}
\newcolumntype{C}[1]{>{\PreserveBackslash\centering}p{#1}}
\newcolumntype{R}[1]{>{\PreserveBackslash\raggedleft}p{#1}}
\newcolumntype{L}[1]{>{\PreserveBackslash\raggedright}p{#1}}

\DeclareMathAlphabet\bpazocal{OMS}{cmsy}{b}{n}
\usetikzlibrary{patterns}

 \usepackage[margin=1in]{geometry}
\usepackage{graphicx}
\usepackage{multirow}
\usepackage{amsmath,amssymb,amsfonts}
\usepackage{amsthm}
\usepackage{mathrsfs}
\usepackage{xcolor}
\usepackage{textcomp}
\usepackage{booktabs}
\usepackage{algorithm}
\usepackage{algorithmicx}
\usepackage{algpseudocode}
\usepackage{listings}
\usepackage{makecell}
\usepackage{authblk}
\usepackage[toc,page]{appendix}
\usepackage{hyperref}

\newtheorem{remark}{Remark}

\title{Leveraging Operator Learning to Accelerate Convergence of the Preconditioned Conjugate Gradient Method}

\author[1,2]{Alena Kopani\v{c}\'akov\'a}
\author[3]{Youngkyu Lee}
\author[3]{George Em Karniadakis}

\affil[1]{Toulouse-INP (ENSEEIHT)/IRIT, 2 Rue Charles Camichel, Toulouse, 31000, France}
\affil[2]{Artificial and Natural Intelligence Toulouse Institute (ANITI), 15 Rue des Lois, Toulouse, 31000, France}
\affil[3]{Division of Applied Mathematics, Brown University, Hope Street 170, Providence, RI 02906, USA}

\date{}

\begin{document}

\maketitle

\begin{abstract}
We propose a new deflation strategy to accelerate the convergence of the preconditioned conjugate gradient~(PCG) method for solving parametric large-scale linear systems of equations.
Unlike traditional deflation techniques that rely on eigenvector approximations or recycled Krylov subspaces, we generate the deflation subspaces using operator learning, specifically the Deep Operator Network~(DeepONet).
To this aim, we introduce two complementary approaches for assembling the deflation operators.
The first approach approximates near-null space vectors of the discrete PDE operator using the basis functions learned by the DeepONet.
The second approach directly leverages solutions predicted by the DeepONet.
To further enhance convergence, we also propose several strategies for prescribing the sparsity pattern of the deflation operator.
A comprehensive set of numerical experiments encompassing steady-state, time-dependent, scalar, and vector-valued problems posed on both structured and unstructured geometries is presented and demonstrates the effectiveness of the proposed DeepONet-based deflated PCG method, as well as its generalization across a wide range of model parameters and problem resolutions.
\end{abstract}


\begin{sloppypar}
\section{Introduction}
Many real-life applications in science and engineering require the numerical solution of parametric elliptic and parabolic partial differential equations~(PDEs).
Solving such systems of equations can be computationally demanding, particularly when the parameters fall within a specific range where high-fidelity solutions are required.
In such cases, it is necessary to solve a sequence of large-scale linear systems, arising after discretization, to a desired level of accuracy.
In this work, we solve these linear systems of equations using the conjugate gradient~(CG) method~\cite{hestenes1952methods}.

The CG method is widely regarded as a robust and efficient iterative method due to its ability to exploit the structure of symmetric positive definite~(SPD) linear systems.
However, its convergence speed tends to deteriorate as the problem size increases~\cite{hestenes1952methods}.
Traditional approaches for improving the convergence of the CG method involve employing suitable preconditioners.
For example, domain-decomposition~\cite{dolean2015introduction} or multilevel~\cite{briggs2000multigrid} preconditioners are commonly used to ensure the algorithmic scalability.
To further enhance the convergence of the preconditioned CG~(PCG) method, we propose hybridizing it with machine-learning~(ML), specifically operator learning approaches.

The idea of using ML to enhance the convergence of iterative methods has recently gained significant attention in the literature.
Researchers have explored both non-intrusive and intrusive hybridization approaches.
Non-intrusive hybridization approaches leverage ML without modifying the internal structure of the algorithm.
This includes for example approaches that use ML to automate parameter selection~\cite{arisaka2023principled,hospedales2021meta,taghibakhshi2023mg,margenberg2022neural, antonietti2021accelerating, katrutsa2020black,grebhahn2016performance}, or to provide a suitable initial guess~\cite{novello2024accelerating,huang2020int,luna2021accelerating,ackmann2020machine}.
Moreover, ML techniques for selecting the most appropriate solution strategy from a given set of methods have been also investigated, see for instance~\cite{bhowmick2006application,motter2015lighthouse}.

In contrast to non-intrusive hybridization approaches, intrusive approaches directly intertwine the algebra with the ML.
On one hand, this requires access to and modification of the source code of the solution strategy under consideration.
On the other hand, it opens a door to the development of novel algorithms.
For example, a hybrid algorithm utilizing a convolutional neural network~(CNN) to approximate the inverse of the discrete Poisson equation -- required to enforce the incompressibility constraint in the operator splitting solver for Navier-Stokes simulations -- has been proposed in~\cite{tompson2017accelerating}.
A different approach was proposed in~\cite{um2020solver}, where an ML model was leveraged to correct errors not captured by the discretized PDE.
Similarly, the authors of~\cite{hsieh2019learning} have introduced a method that modifies the updates of an iterative solver using a deep neural network (DNN).
Another approach, which has been proposed in~\cite{doncevic2022recursively}, focuses on meta-learning the superstructure of numerical algorithms through recursively recurrent neural networks~(RNNs).

A significant focus has also been placed on using hybridization approaches to improve the convergence properties of Krylov methods.
This includes, for example, approaches for enhancing the search directions generated by the CG algorithm~\cite{kaneda2022deep}. 
A large number of methods have been also proposed to improve the quality of the preconditioners. 
For example, the approaches for predicting optimal sparsity patterns of block-Jacobi and incomplete LU~(ILU) preconditioners have been developed in~\cite{gotz2018machine,stanaityte2020ilu}. 
In~\cite{ruelmann2018prospects}, the authors explored learning sparse approximate inverse~(SPAI) preconditioners, while~\cite{li2023learning} focuses on learning approximate matrix factorizations. 
In~\cite{dimola2025numerical}, a neural preconditioner was proposed for mixed-dimensional PDEs, which utilizes a network to approximate the inverse of the system matrix, enabling faster convergence of a Krylov solver.
In~\cite{nieto2024graph}, graph neural networks (GNNs) were employed to construct data-driven preconditioners that adapt to the structure of the linear system, in turn improving the efficiency of the GMRES algorithm.
Similarly, the authors of~\cite{giraud2025neural} present an U-Net based preconditioner, trained in an unsupervised manner, in order to approximate the inverse of the discretized Helmholtz operator.
In the context of domain-decomposition preconditioners, approaches for replacing the discretization and solution process of the subproblems were proposed in~\cite{li2020deep,li2019d3m,lee2025nonoverlapping}. 
Moreover, a significant focus has been given to enhancing the construction of the coarse spaces, see for example~\cite{heinlein2019machine, klawonn2022learning, chung2021learning,heinlein2021combining,ciaramella2022spectral}. 
In the context of multilevel/multigrid preconditioners, the ML has been employed to enhance the design of transfer operators~\cite{weymouth2022data, luz2020learning, greenfeld2019learning, taghibakhshi2021optimization, wang2023learning}, as well as smoothers~\cite{huang2022learning, chen2022meta,kopanivcakova2024deeponet}.  
Furthermore, several approaches that take advantage of spectral bias in order to design effective coarse space solvers have been proposed, see for example~\cite{kopanivcakova2024deeponet,cui2022fourier,azulay2022multigrid,lerer2023multigrid,zhang2022hybrid,kahana2022geometry,lee2025fast}.

In this work, we propose to enhance the performance of the PCG method by incorporating ML-based deflation strategy, giving rise to DeepONet-based deflated PCG~(DPCG) method. 
The key idea behind deflation~\cite{saad2000deflated,frank2001construction} is to accelerate convergence by eliminating components of the solution associated with unfavorable eigenvalues of the preconditioned system matrix.
To this aim, various strategies for constructing deflation operators have been developed in the literature, such as (approximate) eigenvectors~\cite{chapman1997deflated}, subdomain-based methods~\cite{frank2001construction}, and recycling approaches~\cite{daas2021recycling,cortes2018pod}.
However, building effective and computationally feasible deflation operators for complex, real-world applications remains an open challenge.
To to achieve this goal, we propose two complementary strategies for generating deflation vectors by taking advantage of the DeepONet~\cite{lu2021learning, goswami2022physics}. 
First, we adapt the trunk-basis~(TB) approach, originally introduced in~\cite{kopanivcakova2024deeponet} for constructing hybrid preconditioners. 
Second, drawing inspiration from classical recycling techniques~\cite{daas2021recycling}, we build the deflation basis using a set of DeepONet-predicted solutions.
Importantly, since in the both cases the DeepONet is used only to construct the deflation operator, the theoretical convergence guarantees of the DPCG method are retained.

To reduce the computational cost of the DeepONet-based DPCG methods, we explore three different strategies for enforcing the structure of the deflation operator by grouping the degrees of freedom~(dofs). 
To this end, we group the dofs by incorporating problem-specific knowledge, by leveraging the structure of the preconditioner, and by applying clustering to the solution predicted by DeepONet.
Through a series of numerical experiments, we demonstrate that the proposed DeepONet-based DPCG method can significantly improve the convergence, robustness, and applicability of the established PCG method for a wide range of problems, including the Darcy equation with jumping coefficients, the heat equation and the linear elasticity.
Moreover, for all benchmark problems, we demonstrate that the proposed DeepONet-based DPCG generalizes well across a wide range of parameters and problem resolutions.

This paper is organized as follows: 
In Section~\ref{sec:numerical_solution}, we review the DPCG method.
In Section~\ref{sec:onet_deflation}, we provide an overview of the DeepONet and propose strategies for constructing DeepONet-based deflation operators.
In Section~\ref{sec:impl}, we describe the benchmark problems used to test and demonstrate the capabilities of the proposed DeepONet-based DCPG method.
Finally, in Section~\ref{sec:num_results}, we demonstrate the numerical performance of the proposed method.
A summary and a discussion of future work are provided in Section~\ref{sec:conclusion}.

 \section{Model Problem and its Numerical Solution}
\label{sec:numerical_solution}
In many engineering applications, the behavior of a system must be investigated with high fidelity under different scenarios, such as variations in material parameters, boundary conditions, or source terms. 
This work, therefore, focuses on designing novel solution strategies for solving a sequence of linear systems of equations arising from the discretization of elliptic parametric PDEs. 
Let ${\boldsymbol{\theta} \in \boldsymbol{\Theta}}$ be a given parameter vector, where $\boldsymbol{\Theta} \subset \R^P$ and $P \geq 1$. 
The high-fidelity discrete system under consideration has the following form: 
\begin{align} 
\Am({\boldsymbol{\theta}}) \uv ({\boldsymbol{\theta}}) = \fv({\boldsymbol{\theta}}), 
\label{eq:lin_system_of_eq} 
\end{align} 
where $\Am({\boldsymbol{\theta}}) \in \R^{n \times n}$ is the SPD matrix and $\fv({\boldsymbol{\theta}}) \in \R^n$ is the vector, which depends affinely on the parameters $\boldsymbol{\theta}$.
The problems of this type might arise, for example, from the discretization of an elliptic second-order PDEs, such as one describing the linear elastic behavior of a material structure. 
In this particular case, $\Am({\boldsymbol{\theta}})$ would represent the stiffness matrix, $\fv({\boldsymbol{\theta}})$  would stand for the force, and $\uv({\boldsymbol{\theta}}) \in \R^n$ would be the vector of sought nodal displacements.

Efficiently solving problems of the type~\eqref{eq:lin_system_of_eq} is also relevant when dealing with linear time-dependent problems. 
For instance, solving a parabolic PDE using an implicit scheme requires solving the following system of equations at each time step: 
\begin{align}
\underbrace{\bigg(\frac{1}{\Delta \tau} \Mm({\boldsymbol{\theta}}) + \Km({\boldsymbol{\theta}})  \bigg)}_{\Am(\boldsymbol{\theta})} \  \uv^{(t)} ({\boldsymbol{\theta}})  = \underbrace{ \bv^{(t)}({\boldsymbol{\theta}})
+ \frac{1}{\Delta \tau} \Mm({\boldsymbol{\theta}})  \uv^{(t-1)} ({\boldsymbol{\theta}})}_{\fv{(\boldsymbol{\theta})}},  \quad t \geq 1,
\end{align}
where $t$ denotes a time-step index, and $\Delta \tau > 0$ is a time-step size. 
Here, $\Km \in \R^{n \times n}$ denotes the stiffness matrix obtained by discretization in space, $\Mm({\boldsymbol{\theta}}) \in \R^{n \times n}$ is the mass matrix and~$\bv^{(n)} \in \R^n$ is the time-dependent force vector. 
Thus, in this case, one is required to solve as many linear systems with different right-hand sides as many time steps are there, for each choice of parameters, which even further amplifies the need for an efficient large-scale solution strategy.

\subsection{Deflated Preconditioned Conjugate Gradient~(DPCG)}
The computational cost of solving parametric problems can be exorbitant, as it requires a solution of many large-scale linear systems of equations. 
In this work, we aim to accelerate the solution of such problems by utilizing the DPCG method, with a DeepONet-based deflation strategy.
This section provides an algorithmic description of the DPCG method, while the details about DeepONet-based deflation will be discussed in Section~\ref{sec:onet_deflation}. 
To simplify our presentation, for the remainder of this work, we omit explicitly stating the dependence on the parameters $\boldsymbol{\theta}$.

Given and initial guess~$\uv^{(0)}$, the CG method seeks for the approximate solution~$\uv^{(i)}$ of~\eqref{eq:lin_system_of_eq} in the Krylov  subspace $\uv^{(0)} + \pazocal{K}_i(\Am, \rv^{(0)})$, defined as $ {\pazocal{K}_i(\Am, \rv^{(0)}) := \text{span} \  \{ \rv^{(0)}, \Am \rv^{(0)} \ldots, \Am^{(i-1)} \rv^{(0)} \} }$. 
Moreover, on each $i$-th iteration, the residual~$\rv^{(i)}$ is required to be orthogonal to a subspace $\pazocal{K}_i (\Am, \uv^{(0)})$, i.e., we have to ensure that $\fv - \Am \uv^{(i)} \perp \pazocal{K}_i (\Am, \uv^{(0)})$. 
The CG algorithm fulfills these two conditions by constructing the approximation~$\uv^{(i)}$ as
\begin{align}
\uv^{(i)} = \uv^{(i-1)} + \alpha^{(i-1)} \pv^{(i-1)}, 
\end{align}
where the  search direction $\pv^{(i-1)}$ is obtained in recursive manner. 
In particular, on each i-th iteration,  $\pv^{(i)}$ is given by a linear combination of~$\rv^{(i)}$ and~$\pv^{(i-1)}$, i.e., 
\begin{align}
\pv^{(i)} =
\begin{cases} 
\rv^{(i)},   \quad &\text{for}  \ i = 0, \\
\rv^{(i)} + \beta^{(i-1)} \pv^{(i-1)},   \quad &\text{otherwise}. \\
\end{cases}
\end{align}
The parameter~$\alpha^{(i-1)}$ is chosen as  ${\alpha^{(i-1)}=\frac{\langle  \rv^{(i-1)},  \rv^{(i-1)} \rangle}{\langle \pv^{(i-1)},  \Am \pv^{(i-1)} \rangle}}$, i.e., such that the residuals~$\rv^{(i)}$ and $\rv^{(i-1)}$ are orthogonal to each other.  
In addition, $\beta^{(i-1)}$ is obtained by enforcing the conjugacy between~$\pv^{(i)}$ and~$ \pv^{(i-1)}$, i.e.,~${\beta^{(i-1)}= \frac{\langle  \rv^{(i)},   \rv^{(i)} \rangle}{\langle  \rv^{(i-1)},  \rv^{(i-1)} \rangle}}$.

The CG algorithm is well-known for its computational efficiency and low memory requirements.
Moreover, after $i$ iterations, the error of the solution approximation~$\uv^{(i)}$ can be bounded from above~\cite{meurant2006lanczos, saad2003iterative} as 
\begin{align}
\| \uv -  \uv^{(i)} \|_{\Am} \leq 2 \| \uv -  \uv^{(0)} \|_{\Am} \Bigg( \frac{\sqrt{\kappa(\Am)} -1}{\sqrt{\kappa(\Am)}+1}    \Bigg)^{i+1},
\end{align}
where $\kappa(\Am) = \frac{\lambda_{max}(\Am)}{\lambda_{min}(\Am)}$ denotes the condition number of $\Am$. 
Thus, the more ill-conditioned $\Am$ is, the larger $\kappa(\Am)$ becomes, which in turn also implies slower convergence of the CG algorithm.
Here, we also point out that $\kappa(\Am)$ only affects the upper bound of the error, i.e., the worst-case convergence rate. 
The actual convergence speed of the algorithm is also influenced by other factors, such as the distribution of the eigenvalues of $\Am$, the right-hand side $\fv$, and rounding errors~\cite{greenbaum1997iterative}.

\subsubsection{Preconditioning}
To improve the convergence of the CG method, we can employ the SPD preconditioner~$\Mm \in \R^{n \times n}$. 
Using the left preconditioning strategy, the linear system~\eqref{eq:lin_system_of_eq} can be transformed into
\begin{align}
\Mm \Am \uv = \Mm \fv, 
\label{eq:prec_system}
\end{align}
which can be solved using the CG method. 
The error bound is in this case given as
\begin{align}
\| \uv -  \uv^{(i)} \|_{\Am} \leq 2 \| \uv -  \uv^{(0)} \|_{\Am} \Bigg( \frac{\sqrt{\kappa(\Mm \Am)} -1}{\sqrt{\kappa(\Mm \Am)}+1}    \Bigg)^{i+1}. 
\end{align}
Ideally, $\Mm \approx \Am^{-1}$ and $\kappa(\Mm \Am) \approx 1$. 
However, obtaining such $\Mm$ in practice is computationally demanding, and therefore cheaper approximations of $\Am^{-1}$ are often used in practise.
For example, one can employ a few iterations of some stationary method, e.g., Jacobi, Gauss-Seidel, or multilevel/domain-decomposition methods; c.f.,~\cite{saad2003iterative}.

\subsubsection{Deflation}
\label{sec:deflation}
Even when the preconditioned system satisfies \( \kappa(\Mm \Am) \ll \kappa(\Am) \), the presence of a few unfavorable (e.g., extremely small or isolated) eigenvalues can significantly degrade the performance of the PCG method.  
Deflation techniques~\cite{saad2003iterative,marek1995deflation,frank2001introduction} address this issue by projecting out the components associated with such eigenvalues, effectively setting them to zero in the spectrum of the preconditioned operator.  
This targeted spectral modification further reduces the effective condition number and can lead to substantial improvements in convergence.

To formally define the deflation procedure, we assume that there are two transfer operators, namely a restriction operator $\Rm \colon \R^n \rightarrow \R^k$ and its adjoint - prolongation operator - ${\Rm^{\top} =\colon \Pm \colon \R^k \rightarrow \R^n}$, that map data from and to a subspace of size $k$, respectively.
Here, we assume that $\Pm$ has the full rank, and that it contains the information about the extreme eigenvalues. 
Moreover, we define the projection operator~$\boldsymbol{\Pi} \in \R^{n \times n}$, an invertible operator~$\Cm \in \R^{k \times k}$ and the matrix~$\Qm \in \R^{n \times n}$ as 
\begin{align}
\boldsymbol{\Pi} := \Id - \Qm \Am, \qquad \Qm := \Pm \Cm^{-1} \Rm, \qquad \Cm := \Rm \Am \Pm.
\label{eq:all_operators}
\end{align}
The deflation~\cite{chapman1997deflated, saad1997analysis} consists of splitting the approximation space into two complementary subspaces.
Thus, we decompose the solution $\uv$ into two parts, i.e., 
\begin{align}
\uv = (\Id - \boldsymbol{\Pi}) \uv + \boldsymbol{\Pi} \uv.
\label{eq:decomposed_sol}
\end{align}
By exploiting the definition of~$\boldsymbol{\Pi}$ given in~\eqref{eq:all_operators}, the first term in~\eqref{eq:decomposed_sol} can be recast as
\begin{align}
(\Id - \boldsymbol{\Pi}) \uv =  \Qm \Am \uv = \Qm \fv.
\label{eq:proj_defl}
\end{align}
In other words, $(\Id - \boldsymbol{\Pi}) \uv$ can be obtained as $\Pm \boldsymbol{\mu} $, where~$\boldsymbol{\mu}$ is the solution of the following reduced linear system of equations 
\begin{align}
\Cm \boldsymbol{\mu} = \Rm \fv.
\end{align}

To recast the second term in~\eqref{eq:decomposed_sol}, we can take advantage of the fact that~$\boldsymbol{\Pi}$ is the projector; thus, it is idempotent.
This allows us to obtain $\boldsymbol{\Pi} \uv$ from a solution~$\hat{\uv} \in \R^n$ of the following deflated system of equations:
\begin{align}
\boldsymbol{\Pi}^{\top} \Am \hat{\uv} = \boldsymbol{\Pi}^{\top}  \fv.
\label{eq:deflated_system}
\end{align}
Note, the problem~\eqref{eq:deflated_system} is indefinite and therefore admits infinitely many solutions. 
However, all of them satisfy~${\boldsymbol{\Pi} \hat{\uv} = \boldsymbol{\Pi} \uv}$. 
We can solve the deflated system~\eqref{eq:deflated_system} using the PCG method.   
Setting~$\hat{\uv}^{(0)}$ such that~$\hat{\rv}^{(0)} \perp \pazocal{R}(\Rm^T)$ ensures that the generated sequence of iterates~$\{\hat{\uv}^{(i)} \}_{i=1}^{i_{max}}$ satisfies~${\hat{\rv}^{(i)} \perp \pazocal{R}(\Rm^T)}$, for all $0 \leq i \leq i_{max}$. 
Moreover, the post-processed sequence~$\uv^{(i)} := \Pm \boldsymbol{\mu} + \boldsymbol{\Pi} \hat{\uv}^{(i)}$ converges to the solution $\uv$ of the original linear system~\eqref{eq:lin_system_of_eq}.

The error bound of the deflated preconditioned conjugate gradient (DPCG) algorithm is given as~\cite{saad2000deflated, vuik1999deflation}
\begin{align}
\| \uv - \uv^{(i)} \|_{\Am} \leq 2 \| \uv - \uv^{(0)} \|_{\Am} \left( \frac{\sqrt{\kappa(\boldsymbol{\Pi}^{\top} \Mm \Am)} - 1}{\sqrt{\kappa(\boldsymbol{\Pi}^{\top} \Mm \Am)} + 1} \right)^{i+1},
\end{align}
where \( \boldsymbol{\Pi}^{\top} \Mm \Am \) denotes the projected preconditioned operator.  
Hence, if the condition number satisfies
\[
\kappa(\boldsymbol{\Pi}^{\top} \Mm \Am) \ll \kappa(\Mm \Am) \ll \kappa(\Am),
\]
then solving the deflated system~\eqref{eq:deflated_system} shall be more efficient than solving the original or the preconditioned linear system.

The DPCG method is summarized in Algorithm~\ref{alg:deflated_pcg}.  
Note, that it reduces to the standard PCG method when \( \Rm \) and \( \Pm \) are null matrices, and to the classical CG method if in addition \( \Mm = \Im \).  
To ensure that the initial residual satisfies the orthogonality condition \( \rv^{(0)} \perp \pazocal{R}(\Rm^\top) \), the initial guess \( \uv^{(0)} \) must be constructed such that \( \Rm \rv^{(0)} = \boldsymbol{0} \).  
To this end, we modify the user-provided initial guess, denoted by \( \uv^{(00)} \in \mathbb{R}^n \), as
\begin{align}
\uv^{(0)} = \uv^{(00)} + \Qm (\fv - \Am \uv^{(00)}).
\end{align}

\begin{algorithm}[t]
  \caption{Deflated Conjugate Gradient (DPCG)}
 \label{alg:deflated_pcg}
  \begin{algorithmic}[1]
\Require {$\Am \in \R^{n \times n}, \fv \in \R^n, \uv^{(00)} \in \R^n, \Mm \in \R^{n \times n}, \Rm \in \R^{k \times n}, \Pm  \in \R^{n \times k}$}
\State $\Am_c^{-1} = (\Rm \Am \Pm)^{-1}$, $\Qm = \Pm \Am_c^{-1} \Rm$ \COMMENTmine{Deflation subspace setup}
\State  $\uv^{(0)} = \uv^{(00)} + \Qm (\fv - \Am \uv^{(00)})$ \COMMENTmine{Projection of user-specified initial guess}
\State $\zv^{(0)} = \Mm \rv^{(0)} $, where $\rv^{(0)}  = \fv - \Am \uv^{(0)}$
\State $\pv^{(0)} = \zv^{(0)} - \Pm \boldsymbol{\mu}^{(0)}$, where $\boldsymbol{\mu}^{(0)}  = \Cm^{-1} (\Rm \Am \zv^{(0)})$
\While{$i=1, 2, \ldots,$ \textbf{until} convergence} 
\State $\alpha^{(i-1)} = \langle \rv^{(i-1)}, \zv^{(i-1)} \rangle/ \langle \pv^{(i-1)},  \Am \pv^{(i-1)} \rangle$
\State $\uv^{(i)} = \uv^{(i-1)} + \alpha^{(i-1)}  \pv^{(i-1)}$
\State $\rv^{(i)} = \rv^{(i-1)} - \alpha^{(i-1)}  \Am \pv^{(i-1)}$
\State $\zv^{(i)} = \Mm \rv^{(i)} $   \COMMENTmine{Preconditioning step}
\State $\boldsymbol{\mu}^{(i)}  = \Cm^{-1} (\Rm \Am \zv^{(i)})$  \COMMENTmine{Deflation step}
\State $\beta^{(i-1)} =  \langle \rv^{(i)}, \zv^{(i)} \rangle/ \langle \rv^{(i-1)},  \zv^{(i-1)} \rangle$
\State$ \pv^{(i)} = \beta^{(i-1)}  \pv^{(i-1)} + \zv^{(i)} - \Pm \boldsymbol{\mu}^{(i)}$
\EndWhile
\State
\Return  $\uv^{(i)}$
\end{algorithmic}
\end{algorithm}

 \section{Deflation Operator via Operator Learning Approaches}
\label{sec:onet_deflation}
The effectiveness of deflation strategies depends on the spectral information captured by the operator~$\Rm$. 
Ideally, the rows of~$\Rm$ would consist of the eigenvectors associated with the eigenvalues at either end of the spectrum of the preconditioned operator~$\Mm \Am$. 
However, obtaining such eigenvectors is computationally expensive, especially for large-scale problems, and therefore approximation strategies are often employed~\cite{nicolaides1987deflation, dostal1988conjugate}.

In the context of linear parametric equations considered in this work, several strategies for efficiently solving the sequence of parametric systems have been proposed in the literature,  see for example~\cite{saad2000deflated,venkovic2021recycling,carlberg2016krylov}. 
The key idea behind these methods is to extract approximate eigenvectors produced by the Krylov method while solving the systems in the sampling sequence and recycle these vectors to augment the Krylov subspace while solving the subsequent systems in the sequence~\cite{parks2006recycling, van1987iterative}.
In this work, we propose an alternative approach, where the transfer operator~$\Rm$ is created by utilizing operator-learning approach, namely DeepONet.

\subsection{DeepONet}
The DeepONet~\cite{lu2021learning} is an operator learning approach, which can be utilized to approximate a mapping between infinite-dimensional function spaces. 
Following~\cite{jin2022mionet}, let $\{\pazocal{Y}^b\}_{b=1}^{N_b}$ and $\pazocal{V}$ be the infinite-dimensional Banach spaces.
Our goal is to learn the map 
\begin{align}
\pazocal{G} \colon \pazocal{Y}^1 \times \cdots \times  \pazocal{Y}^{N_b} \rightarrow \pazocal{V}. 
\label{eq:DON_multi_dim}
\end{align}
In the context of the parametric linear systems arising from the discretization of PDEs, considered in this work, this can involve various scenarios. 
For example, we can learn a map from a parametrized right-hand side, material parameters or/and boundary conditions to the solution of the underlying PDE.

The DeepONet approximates~$\pazocal{G}$ by utilizing two types of sub-networks, the output of which is combined via inner product operation. 
Branch networks $\{B^b\}_{b=1}^{N_b}$, where $B^b: \R^{ny_b} \rightarrow \R^p$, are used to encode the input functions associated with the PDE parametrization. 
In practice, input to each branch network is represented by the finite-dimensional approximation of the infinite-dimensional input function $y^{b} \in \pazocal{Y}^b$.
In this work, we approximate each $y^b$ in a finite-dimensional space $\pazocal{Y}_h^b \subset \pazocal{Y}$ by evaluating $y^b$ at  $ny_b$ points $\{\qv_j^b\}^{ny_b}_{j=1}$, called sensor locations, which gives rise to a finite-dimensional vector~$\yv^b \in \R^{ny_b}$. 
Note that each input function~$y^{b}$ can be discretized by using a different set of sensor locations~$\{\qv_j^b\}^{ny_b}_{j=1}$.

Trunk network ${T \colon \R^{d} \rightarrow \R^p}$ is used to encode the computational domain.
Thus, the input is a set of coordinates $\boldsymbol{\xi} \in \Omega$, where $\Omega$ denotes the computational domain. 
Combining the output of branch and trunk networks, we can now approximate the nonlinear operator~$\pazocal{G}$ given by~\eqref{eq:DON_multi_dim} as
\begin{align}
\pazocal{G}(y^1, \ldots, y^{nf})(\boldsymbol{\xi}) \approx 
\sum^p_{q=1} \underbrace{B_q^1(\yv^1)  \times \cdots \times B_q^{N_b}(\yv^{N_b})}_{\text{coefficients}} \  \times \underbrace{T_q(\boldsymbol{\xi})}_{\text{basis functions}},
\label{eq:multiinput_DeepONet}
\end{align}
where $\boldsymbol{\xi} \in \Omega$. 
Here, the symbol $B_q^b$ denotes the q-th element of the output of the $b$-th branch network,  while the symbol $T_q$ denotes the q-th output of the trunk network. 
As we can see, the trunk network provides $p$ basis functions, which are then linearly combined with coefficients provided by the branch networks.

\begin{figure}
\begin{minipage}{0.48\linewidth}
\rotatebox{-90}{
\scalebox{0.8}{
\includegraphics{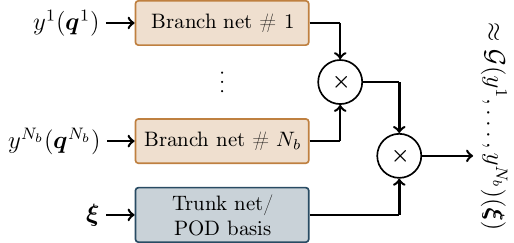}}}
\end{minipage}
\begin{minipage}{0.48\linewidth}
\scalebox{0.72}{
\includegraphics{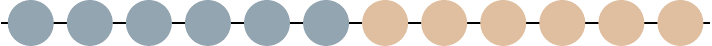}
}
\vspace{0.5cm}

\scalebox{0.72}{
\includegraphics{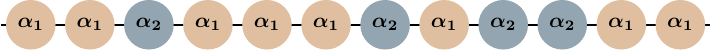}
}
\end{minipage}
\caption{Left: An example of the multi-input DeepONet~\cite{lu2022comprehensive}.
Right: Examples of two groups, illustrated by different colors.
The groups are generated by partitioning the computational domain based on the devision into subdomains (top) or by grouping dofs based on different material properties, denoted by $\alpha_1$ and $\alpha_2$ (bottom).}
\label{fig:DeepONets}
\end{figure}

Figure~\ref{fig:DeepONets} illustrates the DeepONet's architecture.  
The designs of the branch and trunk networks are flexible and can be adapted based on the nature of the input data.  
For example, the branch network may be chosen as a fully connected neural network when the input is a collection of scalar values or as a convolutional neural network when the input consists of multidimensional functions discretized on a structured grid.
Since the coordinate vector \( \boldsymbol{\xi} \in \mathbb{R}^d \) is typically low-dimensional, fully connected architectures are commonly used for the trunk network.  
Alternatively, the trunk network can be replaced with precomputed basis functions, e.g., by applying proper orthogonal decomposition~(POD) to the training data~\cite{lu2022comprehensive} or by parameterizing the integral kernel in Fourier space~\cite{li2020fourier}.

In order to find optimal parameters of the DeepONet, we construct the dataset ${\pazocal{D}=\{(\yv_j^1, \yv_j^2, \ldots, \yv_j^{b},  \bar{\boldsymbol{\xi}}_j,  \uv_j) \}_{j=1}^{N_s}}$.
Here, each $j$-th sample takes into account the discrete input functions $\{ \yv^{b} \}_{b=1}^{N_b}$, and a set of nodal points~${\bar{\boldsymbol{\xi}}_j = [\boldsymbol{\xi}_{j, 1}, \ldots, \boldsymbol{\xi}_{j, n_{\text{don}}}]^{\top} \in \R^{n_{\text{don}} \times d}}$. 
Furthermore, a target solution $\uv_j \in \R^{n_{\text{don}}}$ is obtained by using a high-fidelity discretization method, representing an approximation of $\pazocal{G}$ at the points given in~$\boldsymbol{ \bar{\xi}}_{j}$.
The training is then performed by minimizing the misfit between the output of the DeepONet and the target solution, i.e., 
\begin{align}
\mathlarger{\min}_{\wv \in \R^{np}} \ \ 	
{\mathlarger{\sum}}_{j=1}^{N_s}  \| \ \tilde{\uv}_j - \uv_j \ \|^2.
\label{eq:training_formula}
\end{align}
Here, the DeepONet solution $\tilde{\uv}_j$ is given as 
\begin{align}
\tilde{\uv}_j := \mathlarger{\sum}^p_{q=1} \Bigg( \prod_{b=1}^{N_b} B^b(\wv; \yv^b_j) \Bigg) \times T_q(\wv; \boldsymbol{\xi}_j), 
\label{eq:deeponet_prediction}
\end{align}
where~$\wv$ denotes all parameters of the DeepONet\footnote{For simplicity, the presentation of the DeepONet architecture avoids an explicit dependence on the parameters~$\wv$.}.

We point out that the DeepONet is trained using a preselected set of coordinate points.
However, during the inference, the solution can be approximated at any point within the computational domain by simply evaluating the output of the trunk network for a different $\boldsymbol{\xi}$.
For instance, DeepONet can be trained using a set of coordinates $\{ \boldsymbol{\xi}^C_{j} \}_{j=1}^{n_{\text{don}}}$, associated with a coarse mesh~$\pazocal{T}^C$, while the inference can be performed using $\{ \boldsymbol{\xi}^F_{j} \}_{j=1}^{n}$, associated with a fine mesh~$\pazocal{T}^F$.
This makes the construction of the dataset and the training process cost-efficient, while ensuring that the DeepONet inference remains independent of the discretization strategy.

\subsubsection{Extension of Vanilla DeepONet to Handle Vector-Valued and Time-Dependent Problems}
The DeepONet can be naturally extended to tackle vector-valued outputs, such as those arising in systems of PDEs, e.g., linear elasticity. 
In such cases, the output of the trunk network \( T \colon \mathbb{R}^d \rightarrow \mathbb{R}^{p \cdot d} \) is reshaped or partitioned into \( d \) sub-vectors, each representing the coefficients of the basis functions for one component of the vector field. 
Similarly, each branch network \( B^b \) is designed to output \( p \cdot d \) values, so that the coefficient-basis product in~\eqref{eq:multiinput_DeepONet} is computed for each component of the output vector. 
This yields a component-wise representation of the solution field, given as
\begin{align}
\pazocal{G}(y^1, \ldots, y^{N_b})(\boldsymbol{\xi}) \approx 
\begin{bmatrix}
\sum_{q=1}^p B^1_{q,1}(\yv^1) \cdots B^{N_b}_{q,1}(\yv^{N_b}) \cdot T_{q,1}(\boldsymbol{\xi}) \\
\vdots \\
\sum_{q=1}^p B^1_{q,d}(\yv^1) \cdots B^{N_b}_{q,d}(\yv^{N_b}) \cdot T_{q,d}(\boldsymbol{\xi})
\end{bmatrix},
\end{align}
where \( B^b_{q,r} \) and \( T_{q,r} \) denote the \( r \)-th component associated with the \( q \)-th mode in the branch and trunk networks, respectively.

For time-dependent problems, such as parabolic or hyperbolic PDEs, the input space must include both spatial and temporal coordinates. 
We therefore augment the trunk input to include time, i.e., 
\[
T \colon \mathbb{R}^{d+1} \ni (\boldsymbol{\xi}, t) \mapsto \mathbb{R}^p,
\]
so that the trunk network learns a spatio-temporal basis. 
The branch networks remain without change.
This enables the DeepONet to approximate a mapping ${ \pazocal{G} \colon \pazocal{Y}^1 \times \cdots \times \pazocal{Y}^{N_b} \rightarrow \pazocal{V} }$, where \( \pazocal{V} \subset L^2(\Omega \times [0,T]; \mathbb{R}^d) \).

\subsection{DeepONet-based Deflation Operators}
In this work, we propose to construct the deflation operator \( \Pm \in \mathbb{R}^{n \times k} \) using the DeepONet.  
Following the methodology introduced in~\cite{kopanivcakova2024deeponet}, this construction involves two main steps: extracting deflation vectors from a pre-trained DeepONet and imposing a block structure on the resulting operator.  
Specifically, we explore two distinct approaches for extracting deflation vectors from the pre-trained DeepONet.
In addition, we also investigate several strategies for grouping dofs, which play a critical role in defining the sparsity pattern of~$\Pm$.

\subsubsection{Extracting the Deflation Vectors from the Pre-trained DeepONet}
The construction of the deflation operator begins with the extraction of suitable vectors from the DeepONet.
To this end, we form the tentative deflation operator
\begin{align}
\widetilde{\Pm} = [\vv_1, \vv_2, \ldots, \vv_k],
\label{eq:tentative_operator}
\end{align}
where the vectors $\vv_{l} \in \R^n$, for all $l=1, \ldots, k$ are obtained using one of the following two approaches.

\begin{figure}
\centering
\scalebox{1.25}{
\includegraphics{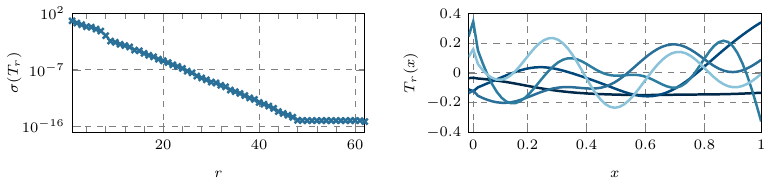}
}
\caption{Left: Visualization of singular values of TB functions extracted from DeepONet ($p=128$), which is trained for a Poisson problem in 1D. 
Right: Visualization of five randomly selected TB functions after performing QR decomposition. }
\label{fig:eigen_fun}
\end{figure}

\begin{itemize}
\item \textbf{Trunk basis~(TB) approach:} 
In the literature, deflation vectors are commonly constructed using (near) null-space vectors~\cite{nicolaides1987deflation}, a strategy often referred to as the Nicolaides approach~(\textbf{NICO}), although it was originally proposed in~\cite{widlund1987additive}.
This approach is particularly convenient in cases where the null-space is explicitly known.
For instance, in the case of the Poisson equation, the null-space contains a constant vector, while in the case of linear elasticity, the null-space vectors consist of rigid body motions related to rotation and translation.
More advanced deflation strategies expand the null-space vectors by incorporating the eigenvectors associated with unfavorable eigenvalues -- specifically, those in the lowest part of the spectrum of~$\Mm \Am$~(see~\cite{frank2001construction}).
In this particular case, $\Pm \Mm \Am$ shares the same eigenvectors as $\Mm \Am$, and the spectrum is given as~\cite{diaz2021accelerating}:
\begin{equation}
\sigma(\boldsymbol{\Pi}^T \Mm \Am) = \{0, \ldots, 0, \lambda_{k+1}, \ldots, \lambda_n \}.
\end{equation}
However, the eigenvectors are typically unknown, and obtaining them is computationally expensive.
As a consequence, approximations are often employed in practice; see, for example,~\cite{burrage1998deflation}.

Our goal is to emulate the behavior of approaches that construct deflation vectors using (approximate) eigenvectors, but without the need for their computationally expensive evaluation.
To achieve this, we extract the vectors~$\{ \vv_{l} \}_{l=1}^k$ using the TB approach from~\cite{kopanivcakova2024deeponet}.
This approach is based on the observation that a large amount of the TB functions is associated with very low singular values, see also Figure~\ref{fig:eigen_fun}. 
In terms of DeepONet approximation properties, these TB functions do not contribute significantly to defining the DeepONet approximation space.
However, they approximate the (near) null-space of the parametric operator well and are therefore well-suited to serve as deflation vectors.

To select the TB associated with the lowest singular values, we can employ the SVD based method, see for example~\cite[Supplement~(Section 2)]{meuris2023machine}. 
However, our numerical experience suggests, that selecting the basis functions at random yields comparable performance in practice.
Therefore, we construct each~$\vv_l$ in~\eqref{eq:tentative_operator} by randomly selecting, without replacement, an index $r < p$, and assembling the vector~$\vv_l$ as
\begin{align}
(\vv_l)_{j} = T_r(\boldsymbol{\xv}_j), \quad \quad \text{for} \quad j = 1, \ldots, n.
\label{eq:matrix_asd}
\end{align}
Here, $T_r(\boldsymbol{\xv}_j) \in \R$ denotes the $r$-th component of the output of the trunk network, evaluated at the coordinate point~$\boldsymbol{\xv}_j \in \Omega$.
In other words, $v_l$ contains $r$-th TB function, evaluated at all nodes of the mesh used for the discretization of the problem at hand~\eqref{eq:lin_system_of_eq}.

\begin{remark}
If the DeepONet is setup with a trunk network consisting of POD basis, the proposed approach can be viewed as a variant of the POD-based deflation outlined in~\cite{cortes2018pod}.
However, in this particular case, the method is not generalizable to problems with varying resolutions.
\end{remark}

\item \textbf{Recycling solutions~(RS) approach:}
Another frequently exploited approach for obtaining deflation vectors is based on so-called recycling strategies~\cite{saad2000deflated,parks2006recycling}.
These strategies exploit the fact that, for a family of linear systems ${ \Am_i \uv_i = \bv_i }$, whose solutions \( \um_i \in \mathbb{R}^n \) span a low-dimensional subspace, any new solution \( \um \) of a related system \( \Am \uv = \fv \) can often be approximated as \( \um \approx \sum_{i=1}^p \alpha_i \um_i \),  
provided that the systems share similar spectral properties or structural dependence.  
This assumption is commonly satisfied for instance if the matrices \( \Am_i \) are obtained through affine parameter sampling.

The solutions \( \{ \um_i \}_{i=1}^p \) can therefore be directly used to construct the deflation matrix \( \widetilde{\Pm} \),  
enabling the projection of error components aligned with slowly converging modes.  
The spectral behavior of the deflated operator is then as follows
\begin{equation}
\sigma(\boldsymbol{\Pi}^\top \Mm \Am) = \{\lambda_1, \ldots, \lambda_{\alpha-1},  0,  \lambda_{\alpha+1}, \ldots, \lambda_{\beta-1},  0,  \lambda_{\beta+1}, \ldots, \lambda_n\},
\end{equation}
where the selected eigenvalues corresponding to the deflation subspace are eliminated from the spectrum.
Here, we highlight the fact that the dominant computational cost of the recycling approaches lies in obtaining suitable \( \{ \um_i \}_{i=1}^p \).

To reduce this computational burden, we propose to predict \( \{ \um_i \}_{i=1}^p \) using a pre-trained DeepONet.  
For \( l = 1 \), the vector \( \vv_1 \) corresponds to the DeepONet predicted solution of the linear system under consideration.  
For \( l = 2, \ldots, k \), predictions are generated by sampling the branch input feature vectors \( \{ \yv^b_l \}_{b=1}^{N_b} \) randomly from the same distribution as used during training.  
Thus, for \( 1 < l \leq k \), each predicted deflation vector \( \vv_l \) is constructed as
\begin{align}
\vv_l = \sum_{k=1}^p \left( \prod_{b=1}^{N_b} B^b(\wv; \yv^b_l) \right) \times T_k(\wv; \xv),
\end{align}
where \( \{ \yv^b_l \}_{b=1}^{N_b} \) denotes the random branch inputs.
In contrast to standard recycling-based deflation strategies, the performance of the RS deflation approach is effective even for the initial systems in the sequence, where no prior solutions are available.

\end{itemize}

\begin{remark}
The proposed DeepONet-based deflation strategies are compatible with both vector-valued and time-dependent problems.  
For vector-valued PDEs, the trunk and branch networks are configured to produce outputs for each component of the solution field.  
The resulting deflation vectors are then constructed by concatenating the component-wise contributions across the domain.
For time-dependent problems, the trunk input is extended to include time, enabling the network to learn spatio-temporal patterns.  
Deflation vectors are then assembled by evaluating the DeepONet either at fixed time instances~(TB approach) or by predicting solutions for a given time-step~(RS approach).  
\end{remark}

\subsubsection{Grouping Degrees of Freedom and Enforcing Block Structure}
Once the tentative transfer operator~$\widetilde{\Pm}$ is constructed using TB or RS approaches, we prescribe the sparsity pattern of~$\Pm$.
This is achieved by grouping the dofs into $S$ disjoint groups.
Each group is associated with an index set~$\pazocal{I}_s$, such that the index set $\pazocal{I}$ of all dofs is given as~$\pazocal{I} = \cup_{s=1}^S {\pazocal{I}}_s$, where~$n := | \pazocal{I} | := \sum_{s=1}^S n_s$, with $n_s := | \pazocal{I}_s |$.

\begin{figure}
\scalebox{0.95}{
\includegraphics{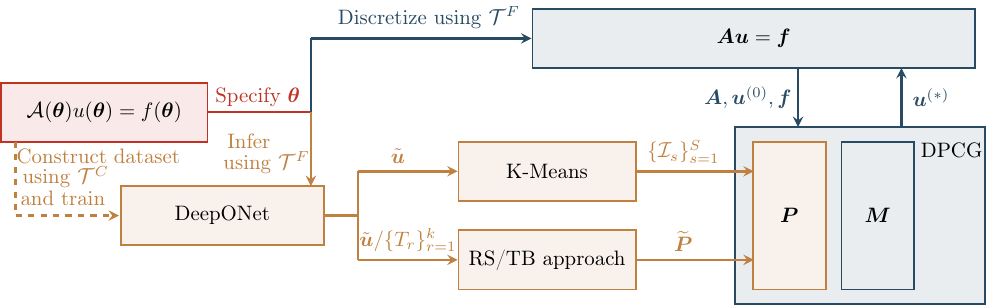}
}
\caption{A sketch of the computational pipeline for DPCG with a DeepONet-induced deflation operator.
The red, blue, and brown colors represent quantities related to the continuous PDE, high-fidelity, and low-fidelity numerical approximations, respectively.
The dashed lines represent the offline DeepONet training stage, while the solid lines are associated with the online stage.
The dataset for training of the DeepONet is constructed using coarse mesh $\pazocal{T}^C$, while at the inference the mesh $\pazocal{T}^F$ of user-desired resolution is utilized.}
\label{fig:pipeline}
\end{figure}

In this work, we explore the following three approaches for constructing~$\{ \pazocal{I}_s \}_{s=1}^{S}$:
\begin{enumerate}
\item \textbf{Problem-specific knowledge:}  
The index sets can be defined based on known features of the problem's structure.  
For example, when a problem exhibits a large jump in coefficients at some location, poor scaling may lead to slow convergence.
In such cases, convergence can often be improved by aligning the group interfaces with the discontinuity, in turn enabling effective deflation.

\item \textbf{Preconditioner structure:} 
The index sets may also be obtained by leveraging the structure of the preconditioner.  
For instance, if the preconditioner is constructed using domain-decomposition approaches, such as the additive Schwarz method~\cite{dolean2015introduction}, then the subdomain partitions used by the preconditioner can be directly utilized to construct $\{ \pazocal{I}_s \}_{s=1}^{S}$.

\item \textbf{Clustering DeepONet predictions:} 
We additionally propose a data-driven approach in which the dofs are grouped by clustering the entries of the solution predicted by the DeepONet.  
Specifically, we use the DeepONet to predict the solution of~\eqref{eq:lin_system_of_eq}, as described in~\eqref{eq:deeponet_prediction}, and pass the resulting vector to a clustering algorithm, namely k-means~\cite{ahmed2020k}.  
Each dof is then assigned to one of the \( S \) groups based on its cluster affiliation.  
In this way, dofs with similar solution behavior are grouped together.  
Consequently, the DeepONet is used not only to construct the deflation vectors, but also to define the block structure and sparsity pattern of the deflation operator.  
An illustration of such DeepONet-based solver pipeline is depicted in Figure~\ref{fig:pipeline}.

\end{enumerate}

Finally, once the index sets have been defined, we perform a block-wise decomposition of the tentative matrix \( \widetilde{\Pm} \) according to these sets.  
As a result, \( \widetilde{\Pm} \) has the following structure:
${\widetilde{\Pm} = 
[
\widetilde{\Pm}_1^{\top},  
\widetilde{\Pm}_2^{\top},  
\cdots, 
\widetilde{\Pm}_{S}^{\top} 
]^{\top},
}$
where each block~$\widetilde{\Pm}_{s} \in \R^{n_s \times k }$ is  associated with the index set~$\pazocal{I}_s$. 
Subsequently, we perform a QR factorization of each block~$\widetilde{\Pm}_s$ to ensure that the deflation basis are linearly independent and well conditioned.
Thus, each block is factorized as \( \widetilde{\Pm}_s = \widetilde{\QQm}_s \widetilde{\Rm}_s \),  
and the orthonormal factor \( \widetilde{\QQm}_s \) is inserted into the global block matrix \( \Pm \in \mathbb{R}^{n \times k \cdot S} \) as follows
\begin{align}
\Pm = 
\begin{bmatrix}
 \widetilde{\QQm}_1 & \boldsymbol{0} & \cdots & \boldsymbol{0} \\
\vdots &  \widetilde{\QQm}_2& \cdots & \boldsymbol{0} \\
\vdots & \vdots& \ddots & \vdots \\
\boldsymbol{0} & \boldsymbol{0}& \cdots &  \widetilde{\QQm}_{S} \\
\end{bmatrix}.
\label{eq:mat_Z}
\end{align}
Thus, each block \( \widetilde{\QQm}_s \), where $1 < s < S$, is inserted into \( \Pm \) such that its row indices align exactly with the corresponding \( \pazocal{I}_s \).

\subsection{Computational Cost of the DeepONet-based DPCG Method}
\label{sec:comp_cost}
The computational cost of the proposed DeepONet-based DPCG method is divided into two stages: offline and online.
In the offline stage, the primary cost is attributed to training the DeepONet model, which is performed only once and amortized afterwords. 
The online stage consists of the initialization and the work performed on each iteration. 
The initialization of DPCG includes the inference through the DeepONet~(RS) or only trunk network~(TS) to generate the deflation basis, construction of the tentative transfer operator \( \widetilde{\Pm} \), QR factorization of its blocks, and assembly of the final deflation matrix \( \Pm \).  
This initialization step incurs a one-time cost.
Each DPCG iteration involves matrix-vector products with \( \Am \), application of the preconditioner \( \Mm \), projections with the deflation operator \( \Pm \), and standard CG vector updates.  
Thus, the per-iteration cost is approximately linear in \( n \).

In practice, selecting a suitable number of deflation vectors \( k \) involves a trade-off between improving the convergence rate of the DPCG method and controlling the associated computational overhead.  
A break-even analysis can be used to estimate the threshold value \( k^\star \) beyond which deflation no longer yields a net computational benefit.  
Assuming that the DPCG iteration count decreases approximately linearly with \( k \), i.e., \( N_{\text{DPCG}}(k) \approx (1 - \theta k) N_{\text{PCG}} \), where \( \theta \) quantifies the effectiveness of each deflation vector, one can derive the value of \( k^\star \) at which the total cost of DPCG matches that of standard PCG.

To estimate \( \theta \) in practice, we can define it as the average relative iteration reduction per deflation vector, i.e., 
\[
\theta \approx \frac{N_{\text{PCG}} - N_{\text{DPCG}}(k)}{k \cdot N_{\text{PCG}}},
\]
where \( N_{\text{PCG}} \) and \( N_{\text{DPCG}}(k) \) denote the number of iterations required to reach a given tolerance using PCG and DPCG with \( k \) deflation vectors, respectively.  
This quantity can be computed empirically by solving the same linear system multiple times with increasing \( k \), and measuring the corresponding convergence rates.
For example, under conservative assumptions (e.g., \( \theta \approx 0.008 \)), deflation remains effective up to \( k \lesssim 30 \), while for moderately effective deflation vectors (e.g., \( \theta \approx 0.01 \)), the benefit extends to \( k \lesssim 50 \).

 \section{Benchmark Problems and Implementation Details}
\label{sec:impl}

\subsection{Benchmark Problems}
This section presents a set of benchmark problems, which we employ for testing and demonstrating the capabilities of the proposed DeepONet-based DPCG algorithm.

\begin{figure}
\hspace{4.0cm}
\includegraphics[scale=0.05]{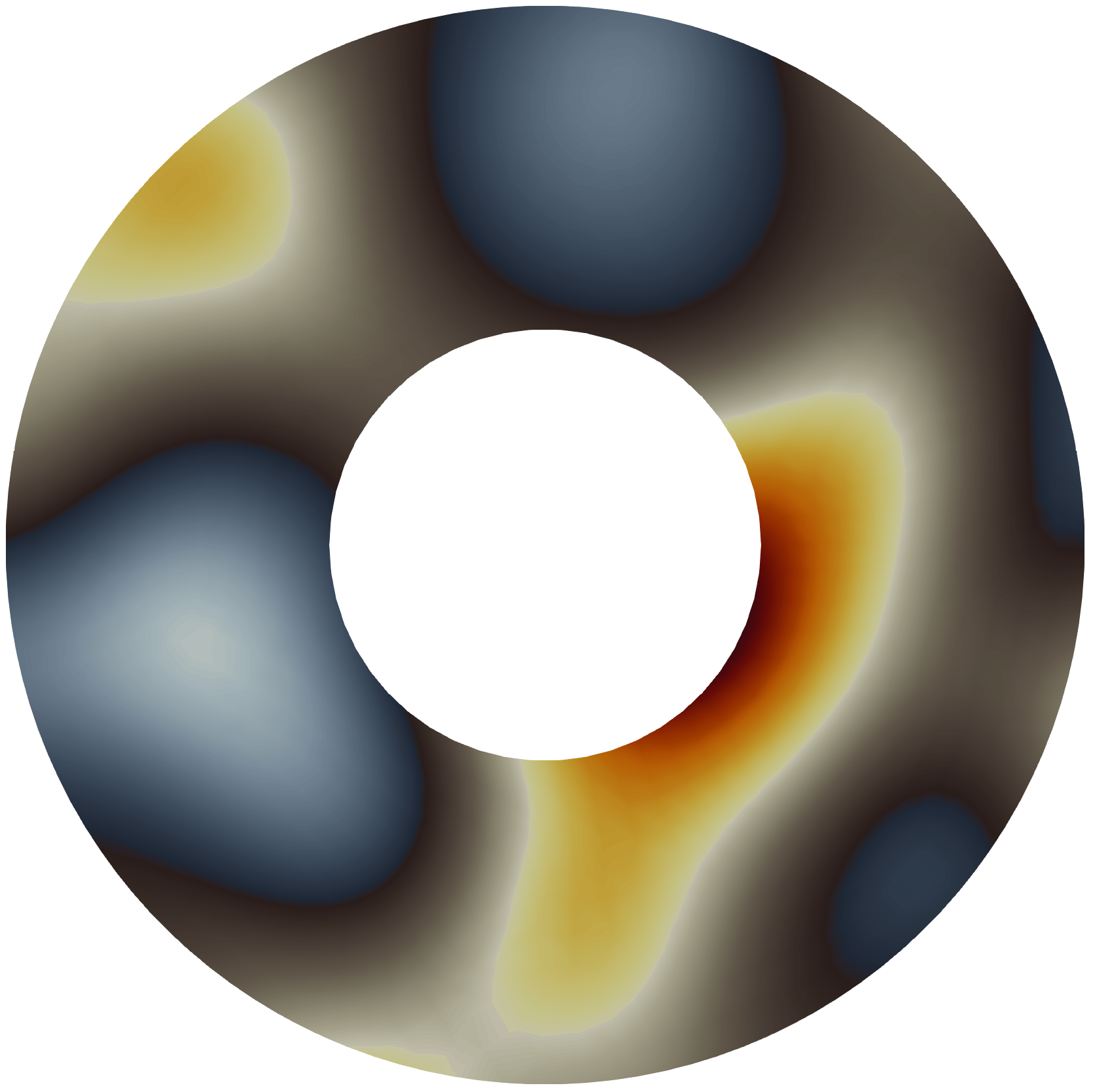}
\hspace{0.5cm}
\includegraphics[scale=0.05]{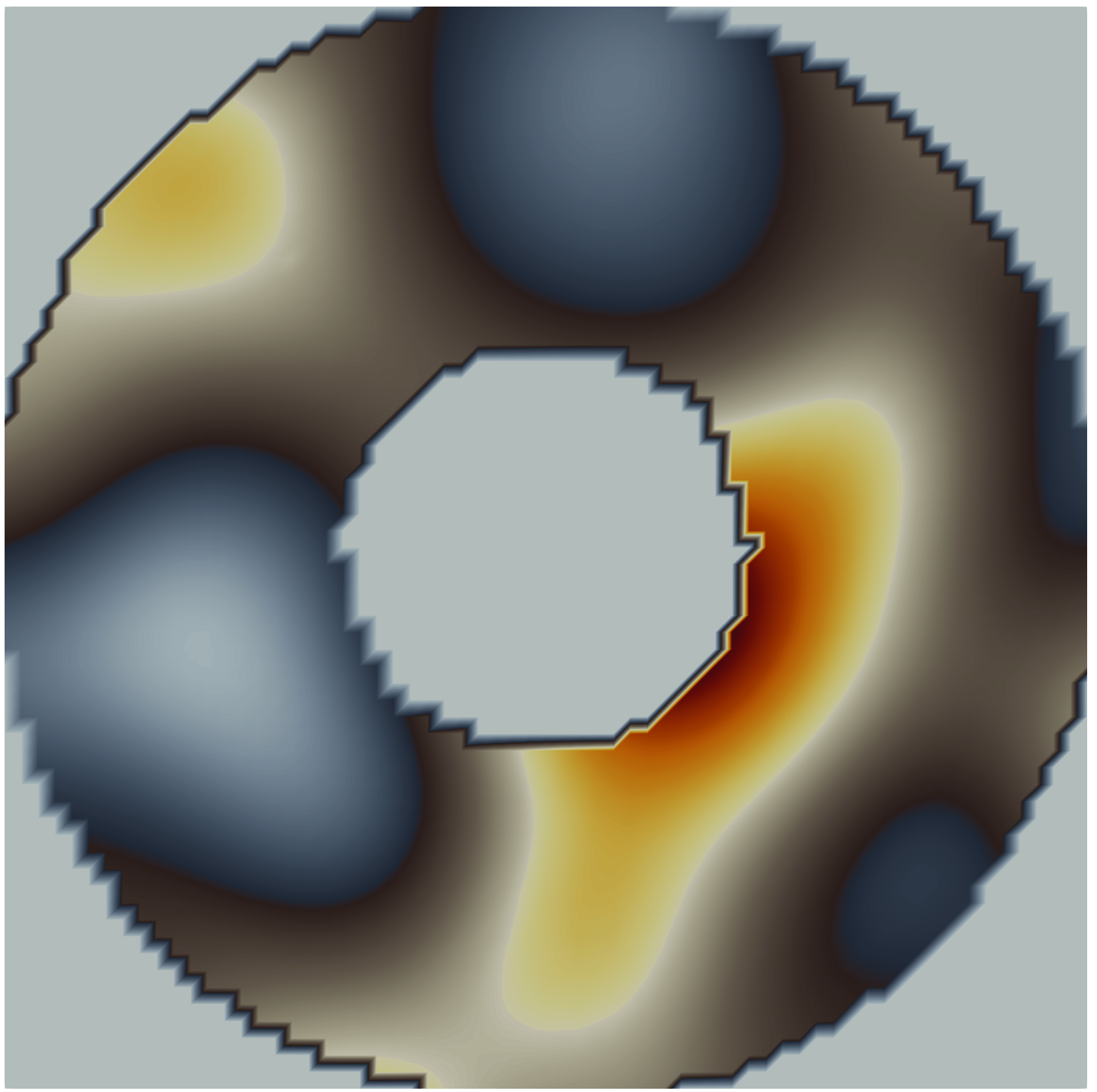}
\hspace{.5cm}
\includegraphics[scale=0.1]{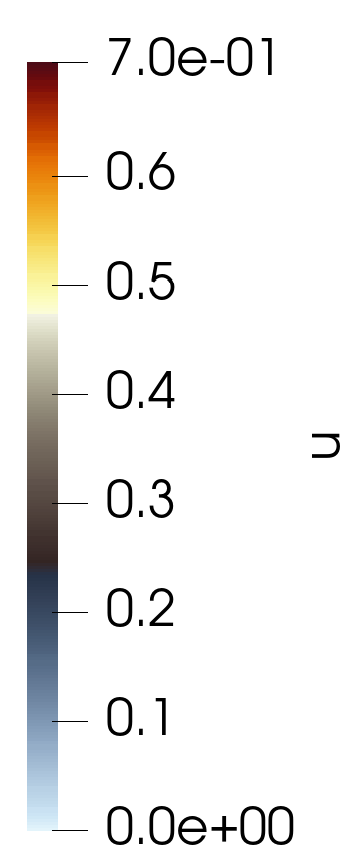}
\caption{Left: An illustration of the spatially varying branch input features ($K$) generated to train DeepONet for Darcy example.
Right: Projection of the branch input features to a bounding box, i.e.,~low-resolution uniform grid, which can be processed by a convolutional branch neural network.}
\label{fig:circle}
\end{figure}

\subsubsection{Darcy Equation}
\begin{figure}
\begin{minipage}{0.23\linewidth}
\scalebox{0.75}{
\includegraphics{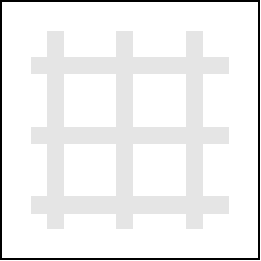}
}
\end{minipage}
\begin{minipage}{0.35\linewidth}{
\includegraphics[scale=0.1275]{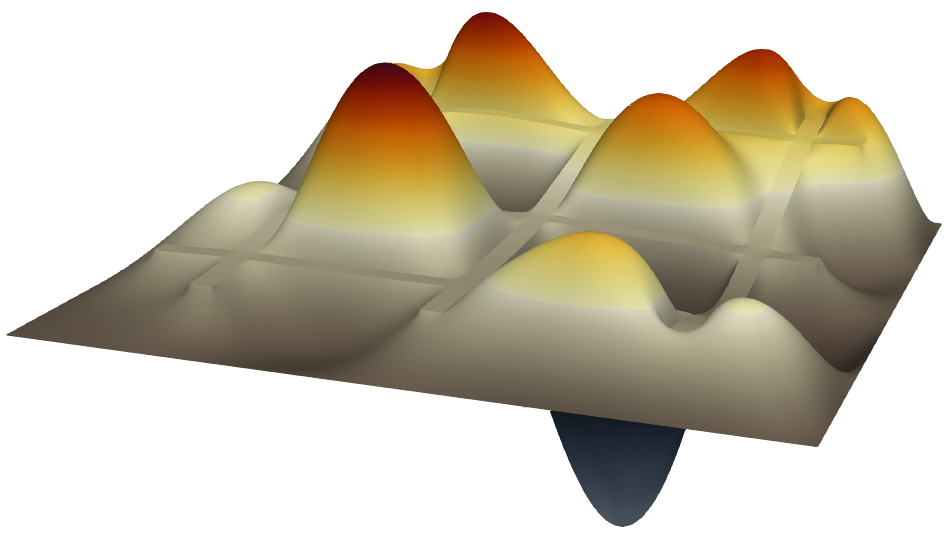}
\includegraphics[scale=0.165]{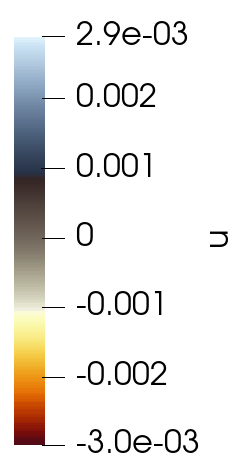}}
\end{minipage}
\hfill
\begin{minipage}{0.3\linewidth}
\includegraphics[scale=0.025]{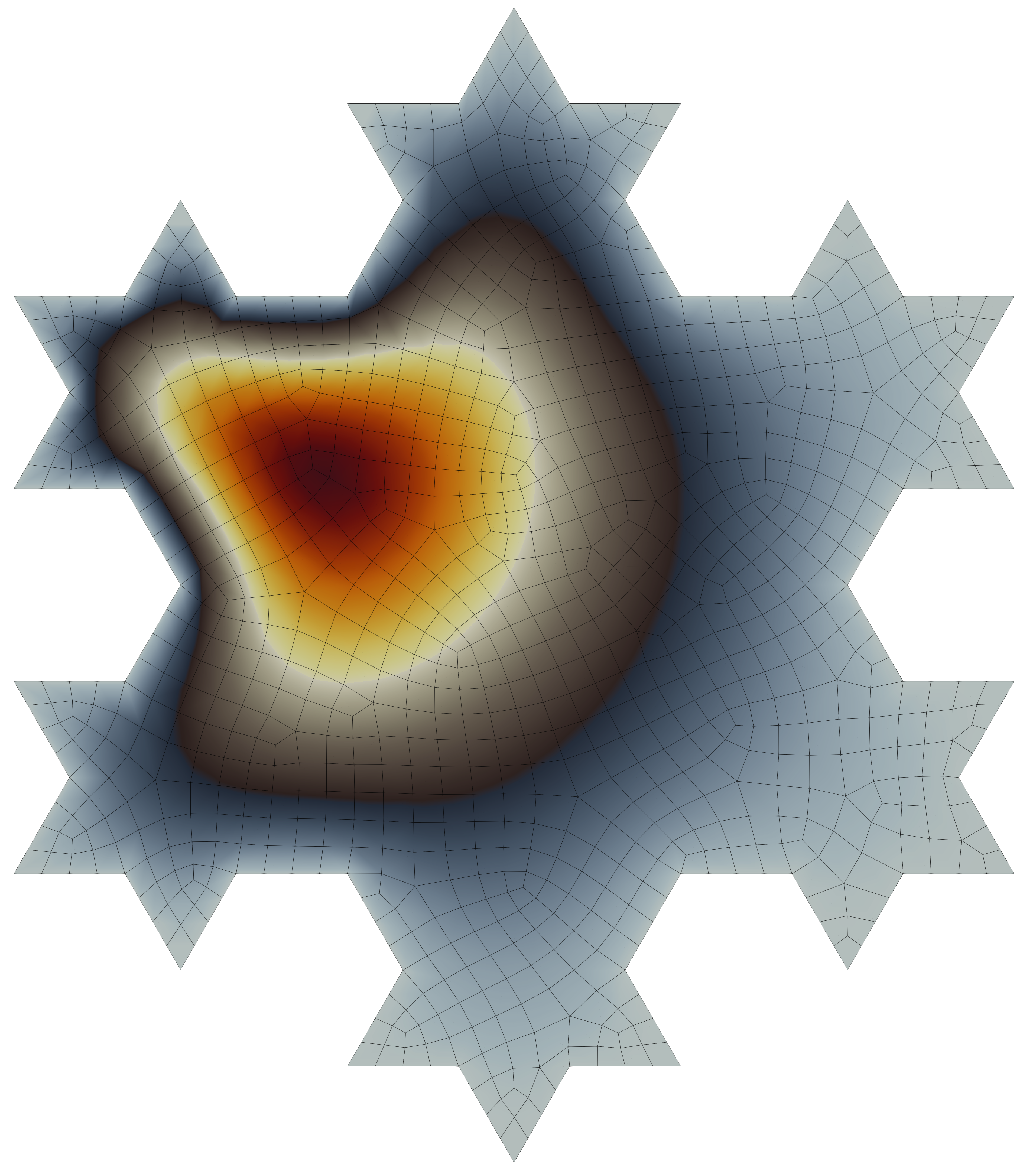}
\includegraphics[scale=0.0455]{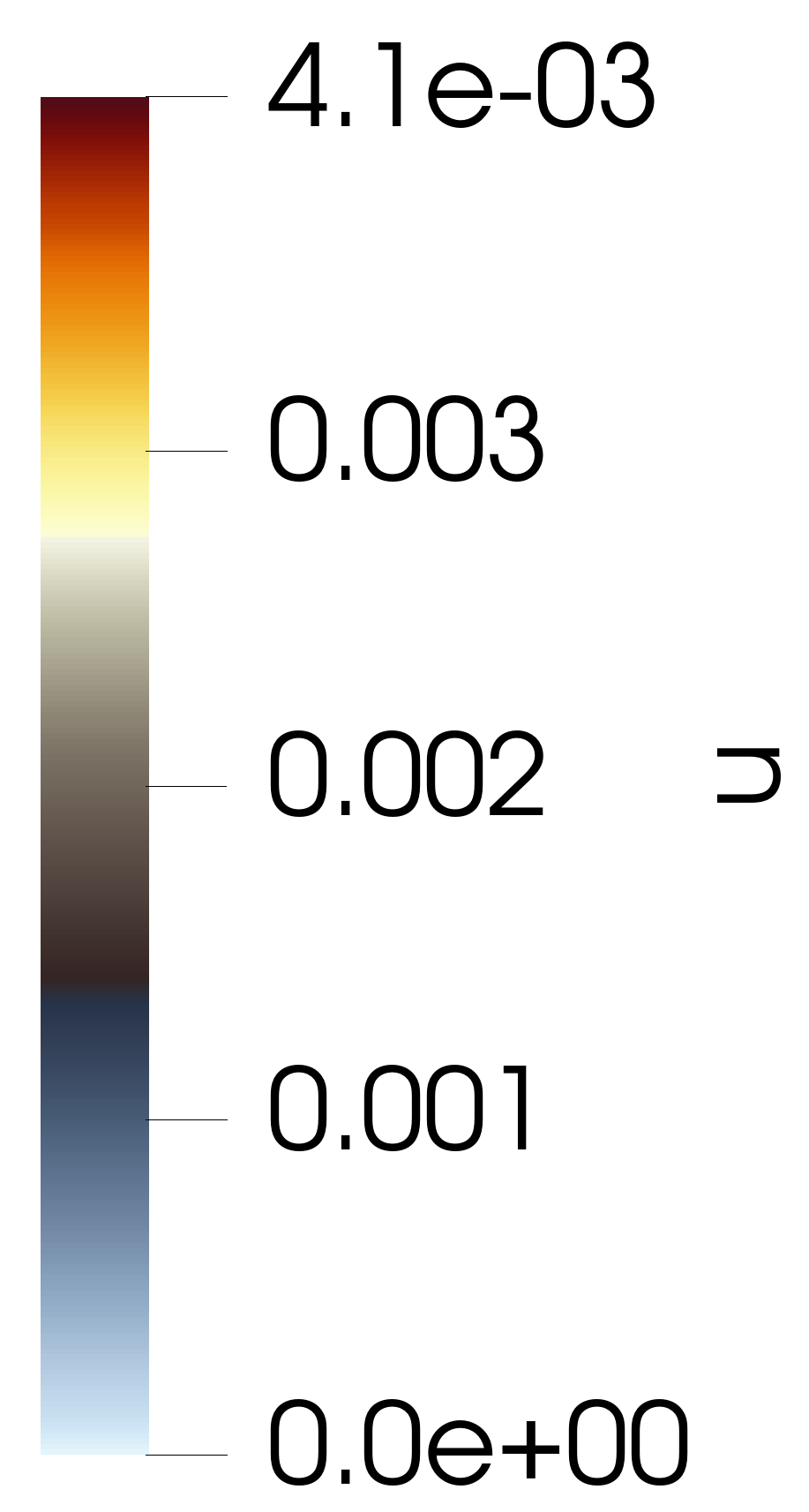}
\end{minipage}
\caption{Left: An illustration of the computational domain with channel patterns used for the Darcy equation test with jumping coefficients. 
The coefficient~$K$ takes on value one in the white region, while in the gray region~$K$ is sampled from~${K \in [1, 10^{5}]}$.
Middle: An example of simulation result for Darcy equation with jumping coefficients. 
Right: An illustration of the initial condition for the snowflake example, depicted on mesh~$\pazocal{T}^1$, which is used to construct the training dataset.}
\label{fig:jump_solutions}
\end{figure}

We start by considering the Darcy equation given as
\begin{equation}
\begin{aligned}
- \nabla \cdot (K(\xv, \boldsymbol{\theta}) \ \nabla u(\xv) ) &= f(\xv, \boldsymbol{\theta}), \qquad  \  \forall \xv \in \Omega, \\
u(\xv) &= 0, \qquad \qquad \quad  \text{on} \ \partial \Omega,
\label{eq:Laplace}
\end{aligned}
\end{equation}
where~$u$ denotes the solution and $f$ stands for the forcing term.
Equation~\eqref{eq:Laplace} is parametrized in terms of the forcing term~$f$ and the diffusion coefficient~$K$.
We consider two different instances of this problem, in particular
\vspace{-0.25cm}
\begin{itemize}
\item \textit{Darcy equation with spatially varying coefficients and forcing term~(Darcy):} 
In this example, we consider an unstructured circular domain~$\Omega$ with radius equal to one, which has a circular hole with a radius of 0.4.
We sample the coefficient~$K$ using Gaussian random fields~(GRFs) with mean $\E[K(\xv, \boldsymbol{\theta})]=0.5$ and  the covariance
\begin{align}
\text{Cov}(K(\xv_1,\boldsymbol{\theta}), K( \xv_2, \boldsymbol{\theta})) = \sigma^2 \exp \bigg(- \frac{\| \xv_1 - \xv_2 \|}{2 \ell^2} \bigg).
\label{eq:sample}
\end{align}
Here, the symbols~$\xv_1$, $\xv_2$ denote the coordinates of two distinct points inside the computational domain~$\Omega$. 
The parameters~$\sigma$ and $\ell$ are chosen as $\sigma = 1.0$ and ${\ell = 0.1}$.
The right-hand side~$f$ is also sampled using GRFs, but with ${\E[K(\xv, \boldsymbol{\theta})]=0.0}$ and the covariance given as in~\eqref{eq:sample}, but with parameters $\sigma = 1.0$ and $\ell = 0.05$. 
The problem is discretized using finite element (FE) method with triangular elements.
An illustration of the geometry, the sampled spatially varying branch input features, and their projection onto the bounding box, which can be processed by the convolutional neural network (branch), is shown in Figure~\ref{fig:circle}.

\item \textit{Darcy equation with jumping coefficients~(JumpDarcy):} 
Next, we consider a problem~\eqref{eq:Laplace} posed at~$\Omega := [0,1]^2$, with fixed right-hand-side ${f(\xv):=\sin(4 \pi \xv_1) \sin(2 \pi \xv_2) \sin(2 \pi \xv_1 \xv_2)}$ and jumping diffusion coefficients. 
The diffusion coefficient~$K$ takes on a value one everywhere, except in grid-like channels, depicted by grey color in Figure~\ref{fig:jump_solutions} (left). 
In channels, the coefficient $K$ takes on a value from $1$ to $10^5$, which we sample as~$\log_{10} K \sim \pazocal{U}[0, 5]$.
The problem is discretized using the FE method, with triangular elements, such that the jumps in the diffusion coefficient are aligned with the edges of the elements.
\end{itemize}

\subsection{Linear Elasticity}
Our next example is associated with linear elasticity in 3D on an unstructured domain $\Omega \subset \R^{3}$, given as
	\begin{equation}
\begin{aligned}
			-\nabla \cdot \sigma\left(\uv,\boldsymbol{\theta}\right) & = \fv(\xv) & \text{ in } \Omega,                        \\
			\uv                                  & = 0                                                                                            & \text{ on } \Gamma,                        \\
			\frac{\partial \uv}{\partial n}      & = \gv(\xv, \boldsymbol{\theta})                                                                                            & \text{ on } \partial\Omega\setminus\Gamma,
		\end{aligned}
        \label{eqn:elasticity}
	\end{equation}
where~$\uv \in \R^3$ denotes the displacement, $\fv$  is the body force, and $\gv \in \R^3$ is the external force.
The stress tensor~$\sigma$ is given as $\sigma\left(\uv,\boldsymbol{\theta}\right) := \lambda(\boldsymbol{\theta}) \mathrm{tr}\left(\varepsilon\left(\uv\right)\right)\Im + 2 \mu(\boldsymbol{\theta}) \varepsilon\left(\uv\right)$, where $\lambda$ and $\mu$ denote the parametrized material parameters. 
The symbol $\Id \in \R^{3 \times 3}$ denotes the identity matrix, and~$\epsilon$ is the linearized strain tensor, given as $\varepsilon (\uv) := \frac{1}{2} ( \nabla \uv + (\nabla \uv)^{T})$.

\begin{figure}
\includegraphics[scale=0.07]{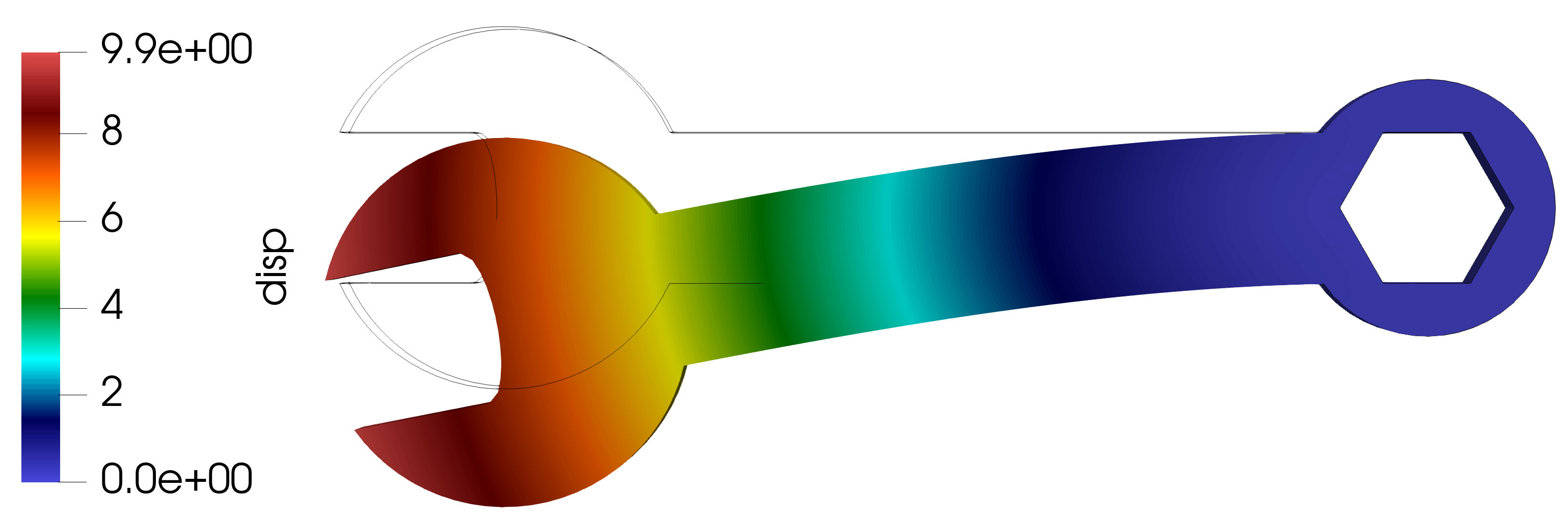}
\hfill
\includegraphics[scale=0.07]{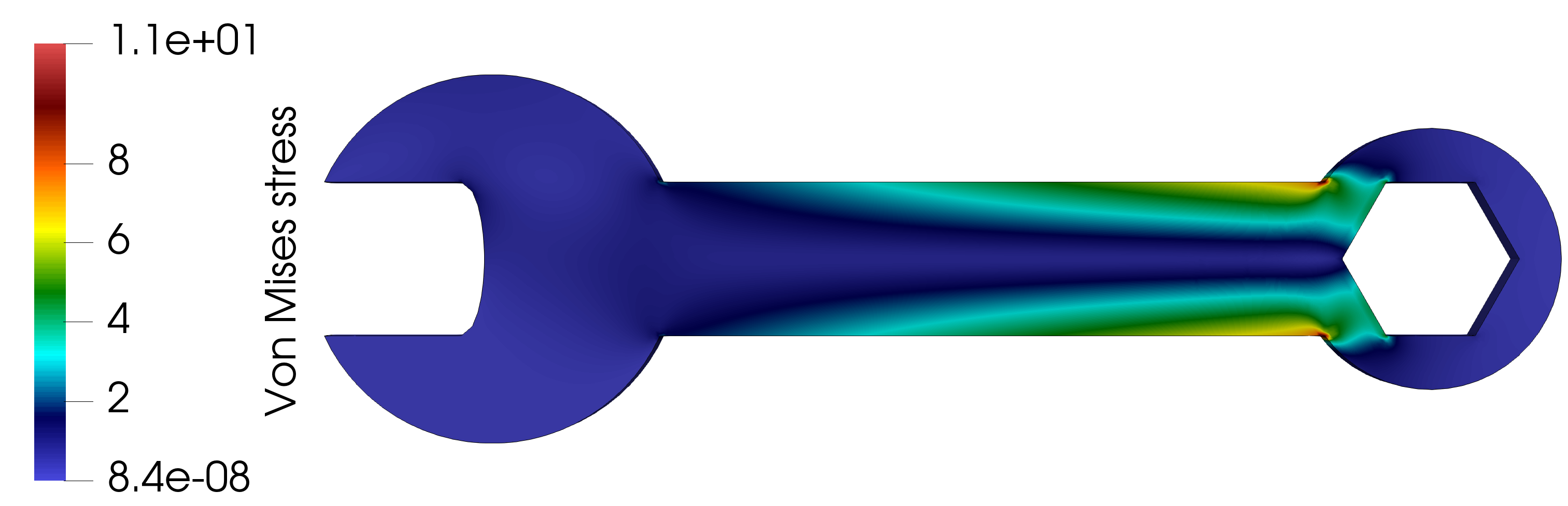}
\caption{A simulation result of a 3D wrench problem under loading with $g_2=0.1$. 
The force is applied on the top of the left jaw of the wrench.
Left: The displacement field. Right: The Von Mises stresses.}
\label{fig:wrench}
\end{figure}

We consider two instances of~\eqref{eqn:elasticity}, i.e., 
\vspace{-0.25cm}
\begin{itemize}
\item \emph{Wrench under varying loading conditions~(Wrench)}\label{sec:wrench}: 
This example considers a 3D wrench with a constant width made out of steel, i.e.,~$\mu = 80$ and~$\lambda = 120$.
The body force $\fv(\xv)$ is as $\fv(\xv) = (0, 0, 0)$, while the parametrized force~${\gv = (0, -g_{2} n_y, 0)}$ is applied on~$\partial \Omega \setminus \Gamma$, which corresponds to the top of the left jaw of the wrench. 
Here, $n_y$ is the surface normal applied in the $y$ direction.
The values of $g_{2}$ are sampled from the distribution $g_{2}\sim\pazocal{U}[-0.05, 0.05]$.
The problem~\eqref{eqn:elasticity} is discretized using the FE method with tetrahedral elements.
An example of the simulation result is shown in Figure~\ref{fig:wrench}.

\item \emph{E-shape with varying material parameters~(E-shape)}\label{sec:eshape}: 
This instance of problem~\eqref{eqn:elasticity} is defined on the 3D E-shaped geometry with $\gv = (0, 0, 0)$, while the body force is set to~$\fv = (0, -0.01, 0)$.
The material parameters~$\mu$ and $\lambda$ are sampled as $\mu\sim\pazocal{U}[25, 80]$ and $\lambda\sim\pazocal{U}[25, 186]$.
As in the previous example, the problem~\eqref{eqn:elasticity} is discretized using tetrahedral FE.
An example of the possible simulation results for a different choice of parameters is shown in Figure~\ref{fig:wrench}.
\end{itemize}

\subsection[Heat]{Heat Equation}
In the end, we consider the heat equation defined on closed and bounded domain $\Omega \in \R^{2}$, i.e., 
\begin{equation}
\begin{aligned}
\frac{\partial u}{\partial t} &= K(\boldsymbol{\theta}) \Delta u, &\text{ in } \Omega \times (0, 1], \\
u &= 0, &\text{ on } \ \partial \Omega \times (0,1], \\
u &= u_{0}, &\text{ at } \ t=0,
\label{eq:Heat}
\end{aligned}
\end{equation}
where $u$ denotes the solution and $K(\boldsymbol{\theta})$ stands for the parametrized thermal diffusivity.
The initial condition is defined as~${u_{0}(x,y) = e^{-5 (x^2+y^2)}}$, and the value of the thermal diffusivity $K$ is sampled as~$K\sim\pazocal{U}[1,2]$.
The problem is discretized in space using the FE method with quadrilateral elements, while in time we use the implicit Euler method with the time step~$\Delta \tau = 0.02$. 

We consider two variations of this problem, which differ from each other by the choice of the computational domain, i.e., 
\vspace{-0.25cm}
\begin{itemize}
\item \textit{Unit square and structured mesh~(Heat):}
For this example, \eqref{eq:Heat} is solved on regular grid ${\Omega = [0, 1]^2}$, discretized using uniform Cartesian mesh.

\item \textit{Snowflake and unstructured mesh~(Snowflake):}
The computational domain for this example is given by the snowflake geometry, generated through three iterations of the Koch snowflake algorithm~\cite{koch1904curve}.
Figure~\ref{fig:jump_solutions} shows the geometry and the unstructured mesh~$\pazocal{T}^1$ employed for generating the dataset.
\end{itemize}

\begin{figure}
\includegraphics[scale=0.0525]{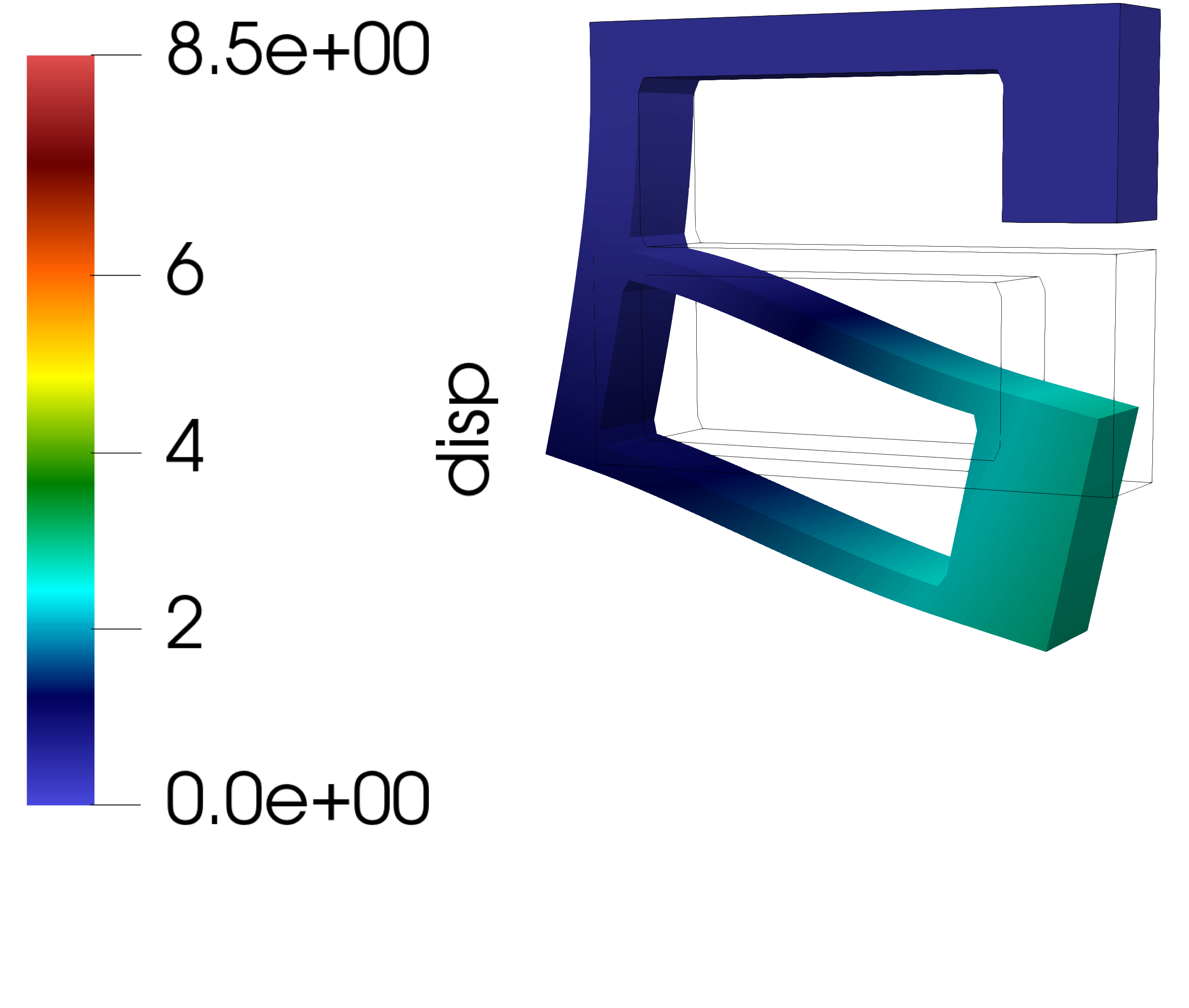}
\hfill
\includegraphics[scale=0.0525]{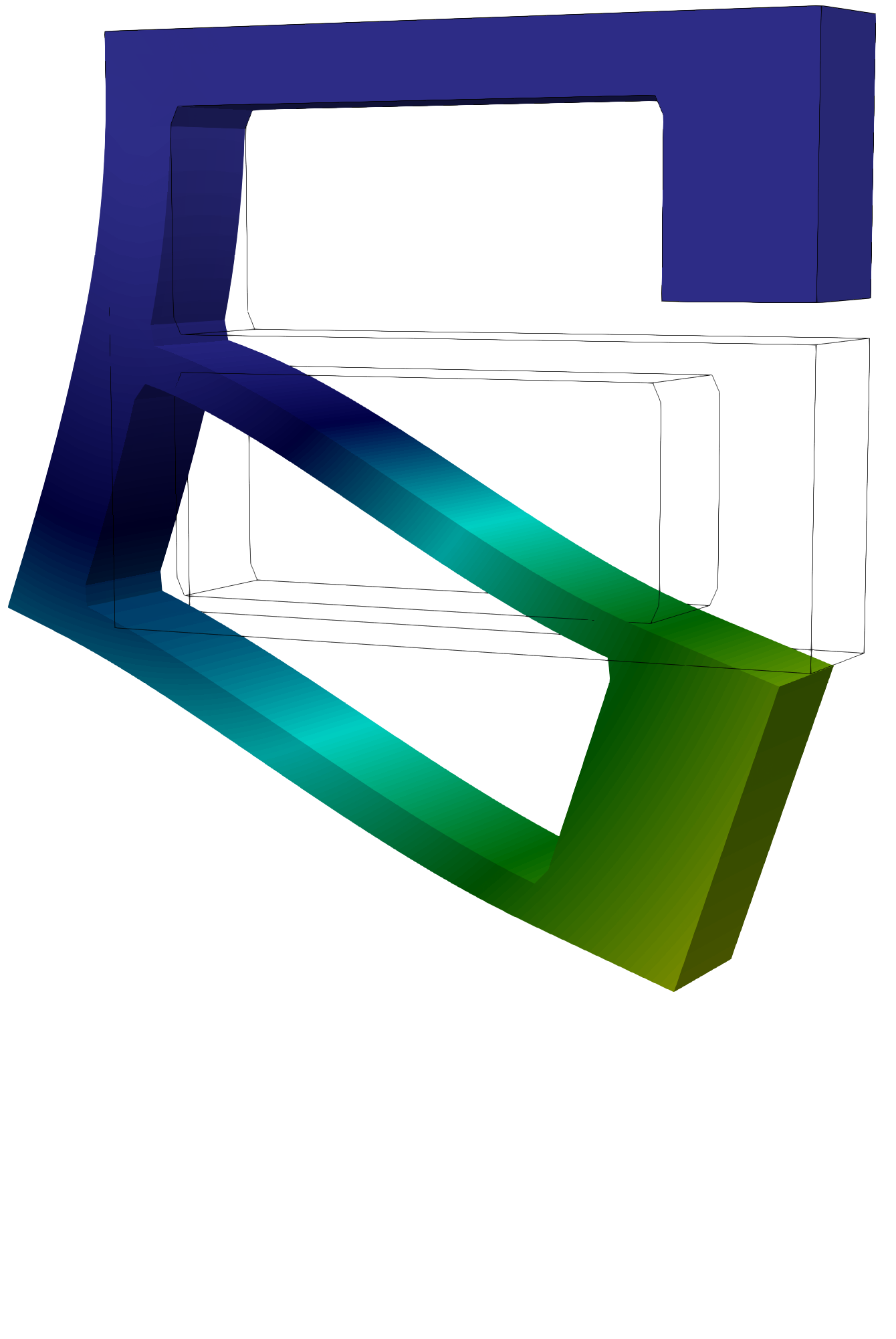}
\hfill
\includegraphics[scale=0.0525]{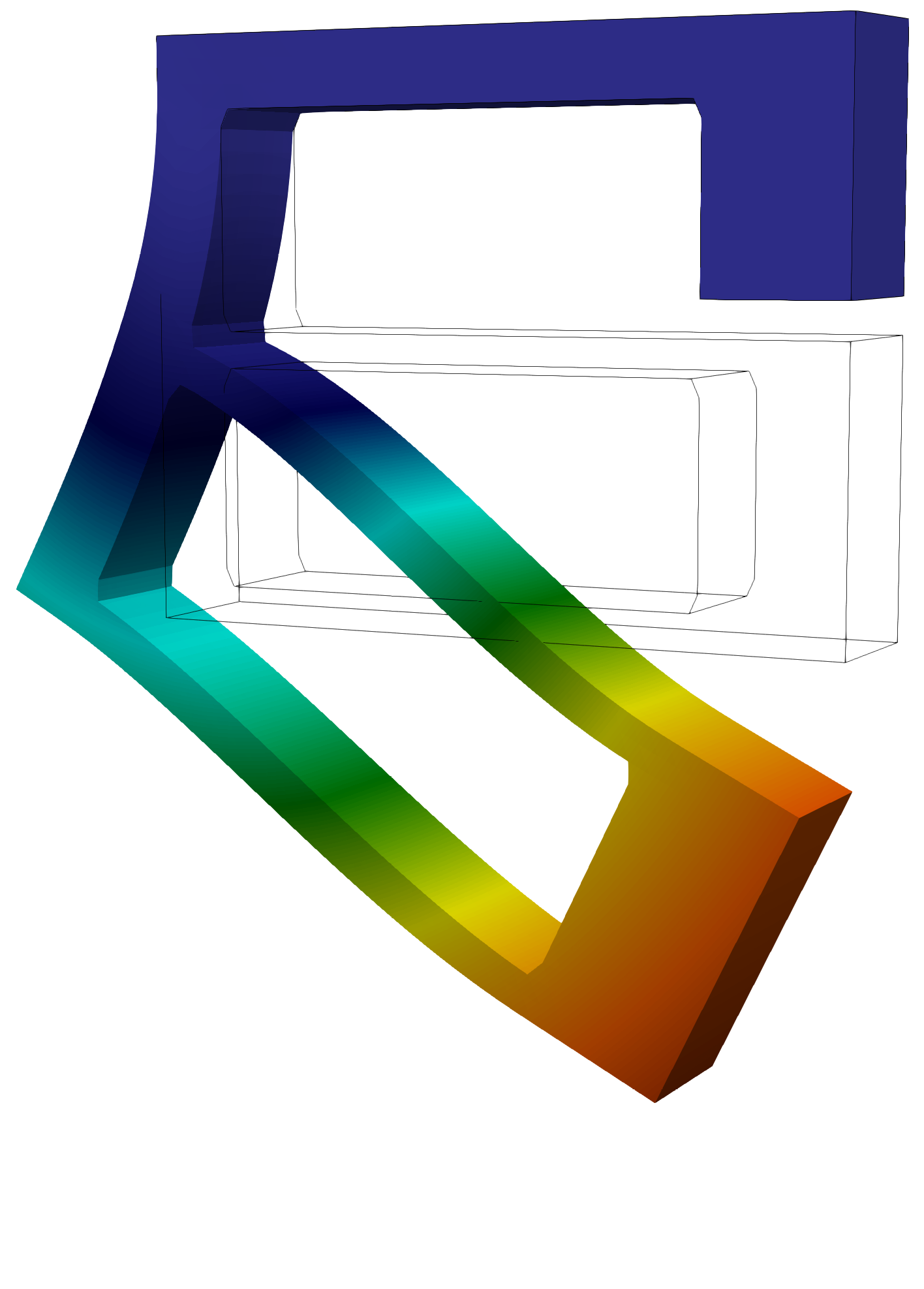}
\hfill
\includegraphics[scale=0.0525]{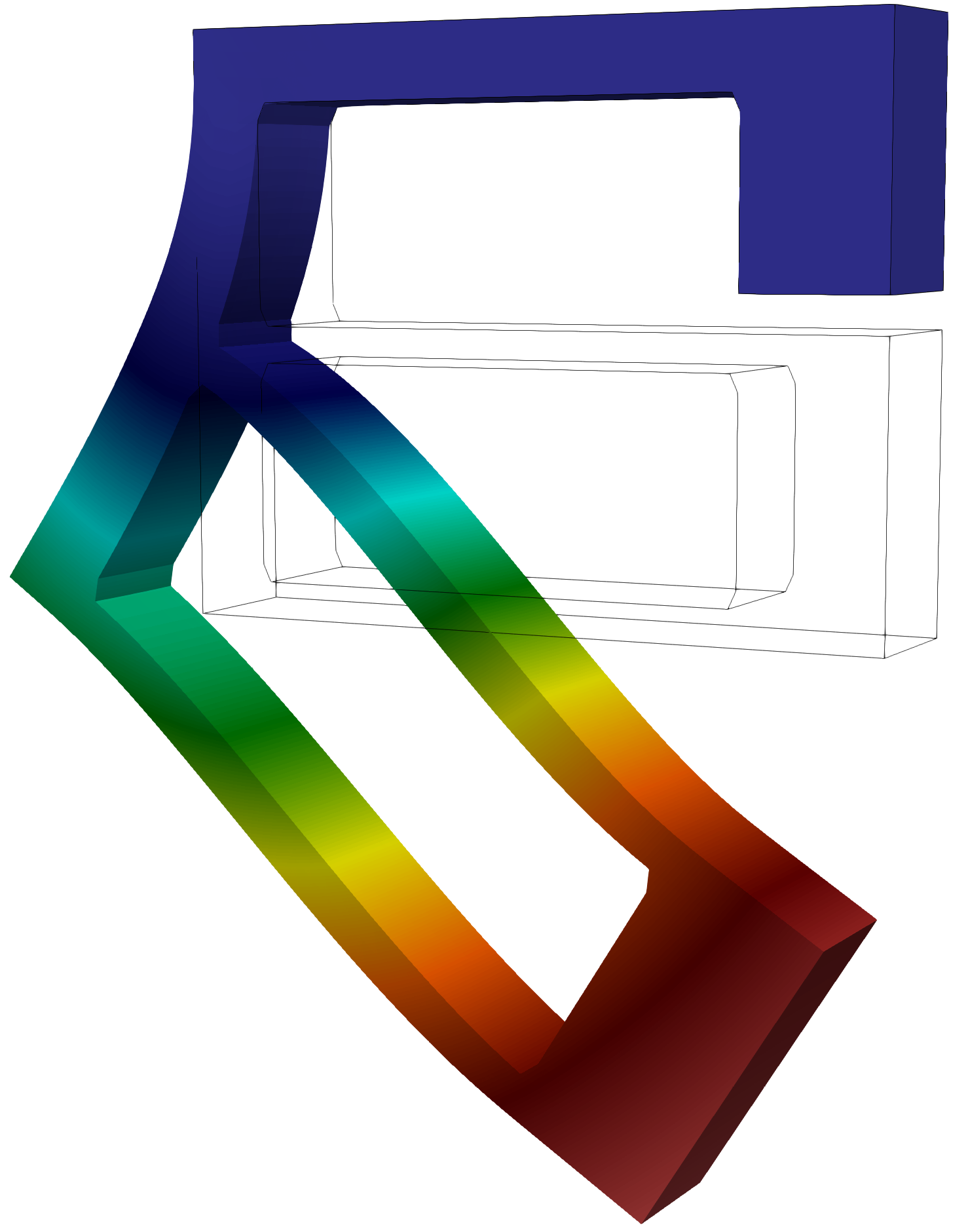}
\caption{An example of a result of a 3D E-shape simulation results for varying material parameters. 
From left to right: $\mu=80, \nu=0.33$;  $\mu=55, \nu=0.27$; $\mu=35, \nu=0.35$; $\mu=30, \nu=0.25$. 
Note, softer materials exhibit larger displacement for the same loading force.}
\end{figure}

\subsection{Implementation Details}
We use the Firedrake library~\cite{rathgeber2016firedrake} to perform the FE discretization of the benchmark problems and to generate the datasets required for training the DeepONets.
The DeepONets are implemented using PyTorch~\cite{paszke2019pytorch} and initialized using the Xavier strategy~\cite{kumar2017weight}.
We train the DeepONets using the Adam optimizer, with a batch size of $1,000$ and a learning rate of $10^{-4}$.
The training process terminates if the validation loss does not improve for $10,000$ consecutive epochs.
Details regarding the network architectures, dataset sizes, and training times are summarized in Appendix~\ref{sec:appB}.
It is worth noting that our numerical experiments are conducted without extensive hyperparameter tuning. 
This suggests that further improvements could be achieved through more elaborate network architecture choices and hyperparameter tuning.
Additionally, the accuracy of the DeepONet and the associated deflation operators could be further improved by employing more advanced training strategies, such as multilevel~\cite{Kopanicakova_2023a, gratton2025recursive} or domain-decomposition methods~\cite{Kopanicakova_2023b,lee2024two}.

The proposed DeepONet-based DPCG methods\footnote{The developed code~\cite{precond_git} will be made publicly available upon acceptance of the manuscript.} are implemented by leveraging the PETSc library~\cite{balay2019petsc}.
The hybridization of the PCG method with DeepONet via deflation is performed using the petsc4py interface.
The numerical experiments were conducted using the Oscar supercomputer at Brown University\footnote{Each computing node is equipped with an AMD EPYC 9554 64-Core Processor (256 GB) and an NVIDIA L40S GPU (48 GB).} and Jean-Zay supercomputer of IDRIS\footnote{Each computing node is equipped with eight V100 GPUs (32 GB) and  24-Core Processor (360 GB).}.

 \section{Results}
\label{sec:num_results}
In this section, we study the numerical performance of the proposed DeepONet-based DPCG method.
The performance is assessed for three different preconditioners: symmetric successive over-relaxation~(SSOR), incomplete Cholesky~(ICC), and the additive Schwarz method~(ASM) with overlap of size one.
Moreover, we consider three different strategies for constructing the deflation matrices: the proposed TB and RS approaches, as well as the standard NICO approach with the deflation vectors constructed as outlined in Appendix~\ref{sec:nico}.

In order to group the dofs, we explore three approaches, i.e.,~computational domain (CD), domain-decomposition (DD) and clustering (CL) approach. 
We explore the CD approach only for the JumpDarcy problem. 
This problem features discontinuous coefficients, which allows us to partition the computational domain based on the two distinct values of the diffusion coefficient~$K$.
Unless specified otherwise, if the ASM preconditioner is employed, we always utilize the DD approach.
Here, the subdomains and the associated index sets $\{ I_s \}_{s=1}^S$ are generated using the Metis partitioner~\cite{karypis1997metis}. 
Notably, the newly introduced CL approach can also be used in conjunction with various preconditioners, including SSOR and ICC.

During all experiments, the DPCG method terminates as soon as one of the following criteria is satisfied:
\begin{align*}
\| \rv^{(i)} \| \leq 10^{-12} \hspace{1cm} \text{or} \hspace{1cm} \frac{\| \rv^{(i)} \|}{\|\rv^{(0)}  \|} \leq 10^{-9}.
\end{align*}
To study the robustness of the DPCG, the performance is evaluated for a wide range of PDE's parameters.
Thus, we always report the number of iterations as an average over $10$ independent runs, i.e., with randomly selected problem parameters and randomly chosen initial guesses.
We also  demonstrate generalization with respect to increasing problem sizes.
For all reported numerical results, the dataset for training DeepONet is constructed using the coarsest mesh~$\pazocal{T}^1$. 
This mesh is then uniformly refined $L-1$ times, giving rise to a hierarchy of $L$ meshes, i.e., $\pazocal{T}^{1}, \ldots, \pazocal{T}^{L}$, which are used for testing the proposed DPCG method.

To demonstrate the convergence of the DPCG, we primarily focus on the algorithmic capabilities and asymptotic convergence.
Thus, we first demonstrate the impact of the possible algorithmic parameters~(number of deflation vectors~($k$), number of the groups~($S$), type of group generator, and the choice of the preconditioner) on the overall convergence of the DPCG method.
Afterward, the comparison with respect to the baseline approach~(PCG without deflation), and DPCG with NICO is provided.

\subsection{Performance with Respect to Increasing Number of Deflation Vectors}
\begin{figure}[t]
\includegraphics[width=\textwidth]{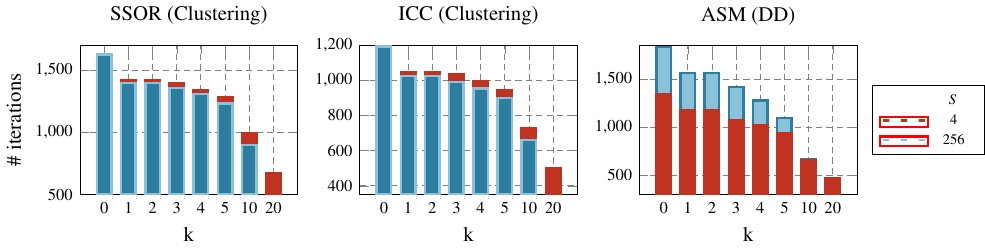}
\caption{The average number of iterations required by the DPCG method with RS deflation to reach convergence. 
	The convergence is monitored with respect to different numbers of deflation vectors~$k$. 
	Note that $k=0$ and $k=1$ denote the baseline (no deflation) and NICO, respectively.
	The experiment is conducted for the Darcy example, discretized using the~$\pazocal{T}^7$ mesh.
	The choices of the preconditioner and the group generator are specified in the title.}
\label{fig:comparison_num_defl2}
\end{figure}

Our first experiment involves studying the convergence behavior with respect to an increasing number of deflation vectors~($k$). Figure~\ref{fig:comparison_num_defl2} presents the results obtained by solving the Darcy problem with SSOR, ICC, and ASM preconditioners and $S \in \{4, 256\}$ groups.
As demonstrated by the results, increasing the number of deflation vectors significantly improves convergence for all DPCG configurations. 
For instance, increasing $k$ from $5$ to $20$ reduces the number of iterations by more than a factor of two, if the ASM preconditioner with $256$ subdomains, i.e.,~$S=256$.
The same convergence behavior was observed across all benchmark problems, which implies that one can enhance the convergence of the DPCG method by simply extracting more deflation vectors from the pretrained DeepONet, albeit at increased iteration cost, see also Section~\ref{sec:comp_cost}.

\subsection{Performance with Respect to Increasing Number of Groups}
\begin{table}
\centering
\begin{tabular}{r|r||r|r||r}
\toprule 
\multirow{3}{*}{\shortstack{Deflation\\type}}   & \multicolumn{1}{r||}{\multirow{3}{*}{S}}  &  \multicolumn{2}{c||}{Darcy} \\
 &   &  \multicolumn{2}{c||}{Preconditioner (Group generator)} \\
  &  & SSOR (CL) & \multicolumn{1}{r||}{ICC (CL) }          \\ \midrule
None                       & 1                    & 1,625.0 & 1,191.6           \\
NICO                       & 4                    & 1,420.4 & 1,043.8           \\
TB             & 4                    & \textbf{1,185.8} & \textbf{869.8}      \\
RS   & 4                    & 1,219.4   & 895.6  \\ \midrule

None                       & 1                    & 1,625.0 & 1,191.6       \\
NICO                       & 16                    & 1,397.4 & 1,028.8 	 \\
TB             & 16                    & \textbf{1,119.4} & \textbf{819.2}  \\
RS   & 16                    & 1,189.2   & 866.8     \\ \midrule

None                       & 1                    & 1,625.0 & 1,191.6       \\
NICO                       & 64                    & 1,398.2   & 1,026.8	 \\
TB             & 64                    & \textbf{1,113.6}   & \textbf{808.2}	 \\
RS   & 64                    & 1,186.4  & 865.6  \\ \midrule

None                       & 1                    & 1,625.0 & 1,191.6          \\
NICO                       & 256                    & 1,397.6  & 1,024.0	\\
TB             & 256                    & \textbf{1,112.4}  & \textbf{808.0} \\
RS   & 256                    & 1,186.4 & 865.6	\\ \bottomrule
\end{tabular}
\begin{tabular}{r|r}
\toprule
  \multicolumn{2}{c}{JumpDarcy} \\
   \multicolumn{2}{c}{Preconditioner (Group generator)} \\
  SSOR (CL) & ICC (CL)           \\ \midrule
  2,481.6 & 1,740.0             \\
  1,966.6 & 1,420.0             \\
  \textbf{1,785.0} & \textbf{1,267.0}     \\
	 1,812.4   & 1,277.4 \\ \midrule

	 2,481.6 & 1,740.0        \\
	 1,952.2 & 1,408.0	 \\
	 \textbf{1,653.2} & \textbf{1,184.6} 	 \\
	 1,673.4   & 1,188.2  \\ \midrule

	 2,481.6 & 1,740.0            \\
	 1,929.4   & 1,389.0	 \\
	 \textbf{1,589.0}   & \textbf{1,142.2}	 \\
	 1,628.8  & 1,165.8	 \\ \midrule

	 2,481.6 & 1,740.0           \\
	 1,912.6  & 1,377.0	\\
	 \textbf{1,563.8}  & \textbf{1,121.4} \\
	 1,617.0 & 1,152.8  \\ \bottomrule
\end{tabular}
\begin{tabular}{r|r||r|r||r|r}
\toprule 
\multirow{3}{*}{\shortstack{Deflation\\type}}     & \multicolumn{1}{r||}{\multirow{3}{*}{S}}  &  \multicolumn{2}{c||}{Heat} \\
&   &  \multicolumn{2}{c||}{Preconditioner (Group generator)} \\
 &  & SSOR (CL) & \multicolumn{1}{r||}{ICC (CL) }          \\ \midrule

None                       & 1                    & 1,719.8 & 1,217.2          \\
NICO                       & 4                    & 1,335.0 & 946.4          \\
TB             & 4                    & \textbf{860.2} & \textbf{617.2}   \\
RS   & 4                    & 874.2   & 624.4   \\ \midrule

None                       & 1                    & 1,719.8 & 1,217.2           \\
NICO                       & 16                    & - & - 	 \\
TB             & 16                    & \textbf{697.4} & \textbf{506.8}  \\
RS   & 16                    & 705.6  & 513.8      \\ \midrule

None                       & 1                    & 1,719.8 & 1,217.2    \\
NICO                       & 64                    & -   & -	 \\
TB             & 64                    & \textbf{641.4}   & \textbf{468.8}	 \\
RS   & 64                    & 649.2  & 473.8  \\ \midrule

None                       & 1                    & 1,719.8 & 1,217.2   \\
NICO                       & 256                    & -  & -	 \\
TB             & 256                    & \textbf{554.0}  & \textbf{403.0} 	 \\
RS   & 256                    & 556.0 & 408.0	\\ \bottomrule
\end{tabular}
\begin{tabular}{r|r}
\toprule 
 \multicolumn{2}{c}{Snowflake} \\
 \multicolumn{2}{c}{Preconditioner (Group generator)} \\
SSOR (CL) & ICC (CL)             \\ \midrule

1,117.8 & 826.2            \\
808.2 & 600.8         \\
583.0 & 436.6   \\
\textbf{578.4}   & \textbf{433.4}   \\ \midrule

1,117.8 & 826.2            \\
 - & - 	 \\
 491.4 & 372.4 \\
 \textbf{481.4}  & \textbf{366.6}     \\ \midrule

 1,117.8 & 826.2    \\
 -   & -	 \\
 457.0   & 347.8	 \\
 \textbf{449.4}  & \textbf{342.0}  \\ \midrule

 1,117.8 & 826.2     \\
 -  & -	 \\
 434.0  & 331.2 	 \\
 \textbf{429.0} & \textbf{328.8}	 \\ \bottomrule
\end{tabular}
\caption{The average number of iterations required by the DPCG method to reach convergence for different types of deflation approaches, preconditioners, numbers of groups (S), and group generators.
For Heat and Snowflake example the average number of the iterations is taken over all time-steps.
The experiment is conducted for the Darcy, JumpDarcy, Heat and SnowFlake examples, discretized using the~$\pazocal{T}^7$, $\pazocal{T}^6$, $\pazocal{T}^5$, and $\pazocal{T}^5$ meshes, respectively.
The number of deflation vectors is chosen to be~$5$ for Darcy and JumpDarcy and $16$ for Heat and Snowflake examples.
A dash (--) is used to indicate that the method failed to converge.}
\label{tab:scalar_valued_groups}
\end{table}

To demonstrate the convergence behavior of the DPCG method with respect to an increasing number of groups, we first consider 
SSOR and ICC preconditioners. 
The grouping of the dofs is performed using the k-means clustering algorithm, while we set $k$ to $5$.
Table~\ref{tab:scalar_valued_groups} presents the results obtained for the Darcy, JumpDarcy, Heat, and Snowflake problems. 
As we can see, increasing the value of $S$ results in slightly improved convergence of the DPCG method. 
For example, for Darcy's problem with the SSOR preconditioner and TB approach, the number of iterations improves only by factor of $1.06$ as $S$ increases from $4$ to $256$. 
The decrease in the number of iterations is more prevalent for time-dependent problems, where the speedup accumulates over multiple time steps. 
However, the reduction in the number of iterations is nevertheless quite moderate, while the iteration's computational cost increases 
due to enlarged size of coarse-space operator~$\Am_c$.
Consequently, when SSOR or ICC preconditioners are used, it is important to keep $S$ relatively low to achieve the tradeoff between the computational cost and the observed speedup.

\begin{figure}[t]
\includegraphics[width=\textwidth]{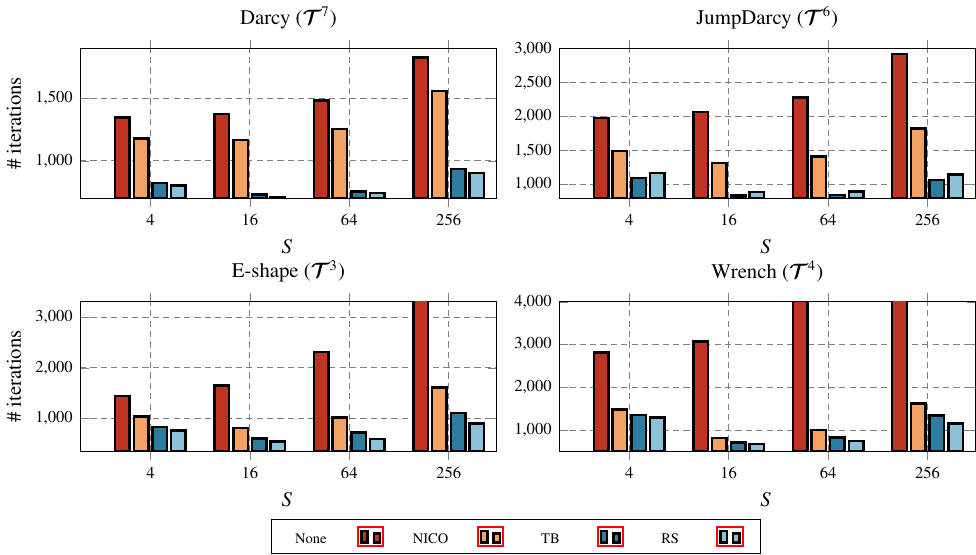}
\caption{The average number of iterations required by the DPCG method to reach convergence with respect to increasing number of groups~(S) associated with number of subdomains of ASM preconditioner. 
We utilized $k=5$~(Darcy and Jump Darcy) and $k=24$~(E-shape and Wrench) deflation vectors for TB and RS approaches, respectively.}
\label{fig:comparison_asm}
\end{figure}

Second, we consider the ASM preconditioner and fix the number of deflation vectors to $k=5$ and $k=24$ for scalar and vector-valued benchmark problems, respectively. 
Here, the groups are determined by the ASM's decomposition into subdomains.
As shown in Figure~\ref{fig:comparison_asm}, the performance of the ASM preconditioner deteriorates significantly when deflation is not applied. 
This is due to the fact that the single-level ASM is not algorithmically scalable, i.e.,~the number of iterations increases as the number of subdomains grows. 
However, incorporating the coarse space via deflation makes the preconditioner scalable for all benchmark problems.
Notably, the TB and RS approaches significantly reduce the number of DPCG iterations compared to the traditional NICO approach. 
This improvement is particularly pronounced for benchmark problems such as JumpDarcy, where the NICO approach fails to produce a suitable coarse space. 
This underscores the practicality of DeepONet-induced coarse spaces for problems where constructing coarse spaces using standard numerical methods is challenging, such as those with jumping coefficients or Helmholtz-like characteristics.
It is worth noting that several robust numerical approaches, such as GenEO coarse space~\cite{spillane2014abstract}, have been developed in the literature to tackle this challenge. 
However, these methods are relatively expensive, particularly during the initialization phase, as they often require solution of an eigenvalue problem for each instance of a parametric PDE.

\subsection{Performance with Respect to Different Group Generators}
Next, we analyze the impact of different group generators on the performance of the DPCG methods for the JumpDarcy example. 
Table~\ref{tab:group_generator} reports the average number of iterations required to reach convergence using various preconditioners and grouping strategies.  
As observed, the best performance for the ASM preconditioner is achieved when combining the TB approach with DD group generator, resulting in a substantial reduction in the iteration count.
For the ICC and SSOR preconditioners, CL grouping outperforms the traditional CD approach.
These findings highlight the fact that the effectiveness of the deflation strategy depends not only on the quality of the deflation vectors but also on the choice of grouping strategy used to impose the sparsity pattern of the deflation operator.

\begin{table}
\centering
%\resizebox{\textwidth}{!}{
\begin{tabular}{r||r|r||r|r||r|r|r}
\toprule
\multirow{3}{*}{\shortstack{Deflation\\type}}	& \multicolumn{7}{c}{Preconditioner type/Group generator} \\ \cmidrule{2-8}
 		& \multicolumn{2}{c||}{SSOR} & \multicolumn{2}{c||}{ICC} & \multicolumn{3}{c}{ASM} \\
    		& \multicolumn{1}{c|}{ CL} &  \multicolumn{1}{c||}{ CD} &  \multicolumn{1}{c|}{ CL}  & \multicolumn{1}{c||}{ CD }	&  \multicolumn{1}{c|}{ CL}	&  \multicolumn{1}{c|}{ DD}  	&  \multicolumn{1}{c}{ CD}    \\ \hline \hline
None  & 2,481.6 		& 2,481.6 		& 1,740.0  		& 1,740.0  		& 2,920.0 		& 2,920.0     		& 2,920.0            \\ 
NICO   & \textbf{1,912.6}  	& 1,944.2 		& \textbf{1,377.0} 	& 1,494.2 			& 	2,278.4 			& \textbf{1,824.4} 	& 2,478.8   \\
TB      & \textbf{1,563.8}  	& 1,944.2 		& \textbf{1,121.4} 	& 1,365.8 			& {1,810.4} 			& \textbf{1,066.0} 	& 1,926.4  \\
RS      & \textbf{1,617.0} 	& {1,848.0}  	& \textbf{1,152.8} 	& {1,316.4} 	          & 1,870.6 			& \textbf{1,146.4} 	& {1,772.2}   \\ \bottomrule
\end{tabular}
%}
\caption{The average number of iterations required by the DPCG method to reach convergence with respect to different type of group generators. 
The number of deflation vectors is chosen to be~$5$ for TB and RS approaches.
The experiment is conducted for the JumpDarcy example, discretized using the~$\pazocal{T}^6$ mesh.}
\label{tab:group_generator}
\end{table}

\subsection{Convergence of DPCG Algorithm with Respect to the Choice of Different Deflation Vectors}
\begin{figure}
	\centering
\includegraphics[width=\textwidth]{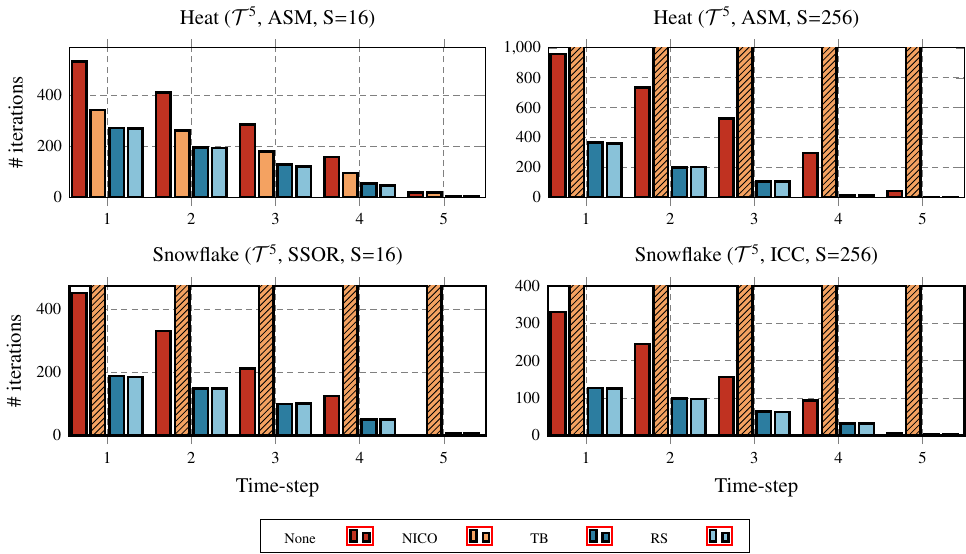}	
\caption{The average number of iterations required by the DPCG method to reach the convergence for different time steps. 
We utilized $k = 16$ deflation vectors for TB and RS approaches.
The mesh, type of preconditioner, and the number of groups are specified in the title.
Note the NICO deflation does not converge for three configurations. 
The bars with a north east lines pattern depict the scenarios where the method did not converge.}
\label{fig:comparison_time_flake}
\end{figure}

Finally, we evaluate the performance of the DPCG algorithm with respect to different choices of deflation vectors.
Table~\ref{tab:asm_summary} summarizes the results for all benchmark problems when using the ASM preconditioner in combination with the DD grouping strategy.  
We observe that both, TB and RS, approaches consistently outperform the standard NICO approach.  
While the performance of TB and RS approaches is generally comparable, the RS approach yields slightly better results across multiple test cases.
However, this behavior is not observed for all preconditioner types.  
In particular, for ICC and SSOR, the results reported in Table~\ref{tab:scalar_valued_groups} show that the TB approach leads to more efficient convergence of the DPCG method.  
Our numerical experiments suggest that this improvement stems from the fact that the TB approach yields a deflation operator that more accurately approximates the spectral components of the preconditioned matrix \( \Mm \Am \).  
This trend is further confirmed in the context of time-dependent problems, where the DPCG method using TB deflation consistently demonstrates improved convergence across all time steps, see also Figure~\ref{fig:comparison_time_flake}.

\begin{table}
\centering
\begin{tabular}{r|r||r|r|r|r|r|r}
\toprule
    \multirow{2}{*}{\shortstack{Deflation\\type}}   &  
    \multirow{2}{*}{\makecell[c]{S}}        	& 
    \multirow{2}{*}{\makecell[c]{E-shape}} 	& 
    \multirow{2}{*}{\makecell[c]{Wrench}}    & 
    \multirow{2}{*}{\makecell[c]{Darcy}}      & 
    \multirow{2}{*}{\shortstack{Jump\\Darcy}}  & 
    \multirow{2}{*}{\makecell[c]{Heat}}       & 
    \multirow{2}{*}{\shortstack{Snow\\Flake}}  \\   
    & & & & & & &  \\ \hline \hline
None                   	& 4        		& 1,444.0      		& 2,817.2   & 1,345.2    & $1,979.6$ & 1,402.4  & 967.8  \\ 
NICO                       	& 4                    & 1,036.0     		& 1,476.6    & 1,179.0    & $1,492.0$ & 895.0 & 631.2\\
TB             		& 4                    & 831.6        		& 1,346.4    & 823.0  & \textbf{1,098.0} & \textbf{617.2}& 490.6 \\
RS   				& 4                    & \textbf{761.0}   	& \textbf{1,295.0}  & \textbf{801.4} & 1,171.4 & 628.4 & \textbf{481.8} \\ \midrule
None                   	& 16                    & 1,647.0        	& 3,074.2    & $1,374.8$ & 2,067.8 & 1,490.4 & 1,072.4 \\ 
NICO                       	& 16                  & 811.4          		& 812.4   & $1,168.2$ & 1,316.2 & 738.2 & 600.6 \\
TB             		& 16                  & 598.4          		& 706.8   & 729.8 & \textbf{840.4} & 519.0 & 442.8 \\
RS   				& 16                  & \textbf{541.8}  	& \textbf{670.6}   & \textbf{709.4}   & 890.4 & \textbf{499.0}   & \textbf{441.6}  \\ \midrule
None                   	& 64                    & 2,310.6     		& $4,015.8$     & $1,482.8$  & 2,280.8 & 1,674.0 & 1,406.8  \\ 
NICO                       	& 64                  & 1,019.0     		& $996.4$     & $1,253.4$ & $1,412.4$ & 806.6 & 763.2 \\
TB             		& 64                  & 716.4       		& $818.6$    & $753.2$ & \textbf{845.2} & 520.2 & 534.4 \\
RS   				& 64                  & \textbf{596.8}  	& $\textbf{742.2}$  & \textbf{741.0} & 898.4 & \textbf{504.0} & \textbf{537.4}  \\ \midrule
None                   	& 256                    & 3,663.6      		& $6,460.0$    & $1,830.4$  & 2,920.0  & 2,559.4 & 2,327 \\ 
NICO                       	& 256                & 1,607.4      		& $1,616.6$    & $1,561.2$ & 1,824.4 & -- & --  \\
TB             		& 256                & 1,108.8      		& $1,333.0$     & $932.6$ & \textbf{1,066.0} & 682.8 & 845.6 \\
RS   				& 256                & \textbf{894.0}  	& \textbf{1,146.8}  & \textbf{899.3} & 1,146.4 & \textbf{668.0 }& \textbf{837.2}  \\ \bottomrule
\end{tabular}
%}
\caption{The average number of iterations required by the DPCG method, preconditioned with ASM, to reach convergence is reported for different types of deflation approaches.
For time-dependent problems, the average is computed over all time steps.
The number of groups is determined by the domain-decomposition of the ASM preconditioner.
For the TB and RS approaches, the number of deflation vectors is chosen to be~$5$ for Darcy ($\pazocal{T}^7$), JumpDarcy ($\pazocal{T}^6$), Heat ($\pazocal{T}^5$), and Snowflake ($\pazocal{T}^5$), and~$24$ for Wrench ($\pazocal{T}^4$) and E-shape ($\pazocal{T}^3$) examples.
A dash (--) is used to indicate that the method failed to converge.}
\label{tab:asm_summary}
\end{table}

 \section{Summary}
\label{sec:conclusion}
We presented a novel framework for accelerating PCG methods using operator learning, with a particular focus on DeepONet-based deflation.
Two strategies for constructing deflation vectors were proposed: (i) the trunk basis~(TB) approach, which exploits the spectral properties of the DeepONet's trunk basis functions, and (ii) the recycled solution~(RS) approach, which utilizes solutions predicted by the DeepONet from randomly sampled input features associated with the PDE parameterization.
To improve the effectiveness of the deflation operator, we investigated three strategies for imposing its block structure: incorporating problem-specific knowledge, leveraging the structure of the preconditioner, and applying clustering to the solution predicted by DeepONet.
Extensive numerical experiments demonstrate that DeepONet-based deflation can substantially reduce the number of PCG iterations across a wide range of problems, including both steady-state and time-dependent scalar- and vector-valued problems, such as parametrized Poisson equations with discontinuous coefficients, linear elasticity, and the heat equation, all defined on various domains and discretized using structured and unstructured meshes.
Notably, the proposed TB and RS approaches consistently outperform standard near null-space-based deflation vectors~(NICO). 
Moreover, the proposed clustering-based grouping of dofs proves particularly effective when no explicit structure in the problem or preconditioner is known a priori.

In the future, we plan to investigate alternative operator learning frameworks, such as Fourier Neural Operators~(FNOs)~\cite{cui2022fourier,li2020fourier} or transformer-based architectures~\cite{pmlr-v235-wu24r,luo2025transolver}, which may yield more expressive and robust deflation bases, especially for multi-physics applications.
We also aim to perform a theoretical analysis of the approximation properties of DeepONet-induced deflation spaces to better understand their effectiveness.
Finally, we envision extending the proposed deflation framework to nonlinear settings with a particular focus on nonlinear Krylov methods.
 
\section*{Acknowledgements}
The work of A.K. benefited from the AI Interdisciplinary Institute ANITI, funded by the France 2030 program under Grant Agreement No. ANR-23-IACL-0002. 
Y.L. was supported in part by Basic Science Research Program through NRF funded by the Ministry of Education (No. RS2023-00247199). 
G.E.K. is supported by the ONR Vannevar Bush Faculty Fellowship. We also acknowledge support from the DARPA DIAL Program, and DOE-MMICS SEA-CROGS DE-SC0023191 award. 
The numerical results were partially obtained using HPC resources from GENCI-IDRIS (Grant No. AD011015766).

\section*{Author contributions}
A.~K.: Conceptualization, methodology, software, writing, reviewing, and supervision. 
Y.~L.: Software, writing, and editing. 
G.~E.~K.: Conceptualization, editing, and supervision.

\begin{appendices}
\section{Numerical Approximation Details}
\label{sec:appB}
In this section, we provide details regarding the high- and low-fidelity numerical approximations for all benchmark problems considered in Section~\ref{sec:impl}. In particular, Table~\ref{table:mesh_hierarchies} summarizes information about the meshes used for FE discretization.
Unless specified otherwise, we use meshes~$\pazocal{T}^1$ to construct the dataset required for training the DeepONets. 
The details of the DeepONet architectures are presented in Table~\ref{tab:architecture}.
For convolutional networks (Conv), the kernel size is always set to three, while the stride is set to two. The feed-forward networks (FFNs) employ standard dense layers consisting of weights and biases.
Table~\ref{table:training_params} provides the required training times for all employed DeepONets and the sizes of the datasets used.

\begin{table}
\centering
\tabcolsep=0.175cm
\begin{tabular}{r||r|r|r|r|r|r|r}
\hline
\textbf{Example}  &  \raisebox{\dimexpr-\height + 1.5ex\relax}{$\bpazocal{T}^1$}	& \raisebox{\dimexpr-\height + 1.5ex\relax}{$\bpazocal{T}^2$}  & \raisebox{\dimexpr-\height + 1.5ex\relax}{$\bpazocal{T}^3$}  & \raisebox{\dimexpr-\height + 1.5ex\relax}{$\bpazocal{T}^4$} & \raisebox{\dimexpr-\height + 1.5ex\relax}{$\bpazocal{T}^5$} & \raisebox{\dimexpr-\height + 1.5ex\relax}
{$\bpazocal{T}^6$} & \raisebox{\dimexpr-\height + 1.5ex\relax}
{$\bpazocal{T}^7$}
\\ \hline  \hline
Darcy 					& $1,106$		& $4,248$		& $16,640$	& $65,856$	& $262,016$ 	& $1,045,248$ & $4,175,360$    \\ \hline
JumpDarcy 				& $2,500$		& $9,801$		& $38,809$	& $154,449$	& $616,225$  & $2,461,761$  &  \\ \hline
Wrench 			& $9,216$		& $57,750$	& $402,516$	& $2,991,048$	& & &  \\ \hline
E-shape 			& $17,388$	& $116,178$	& $843,084$	& & & &  \\ \hline
Heat					& $2,500$		& $9,801$		& $38,809$	& $154,449$	& $616,225$  & & \\ \hline
Snowflake					& $1,007$		& $3,833$		& $14,945$	& $59,009$	& $234,497$  & & \\ \hline
\end{tabular}
\caption{Summary of the number of spatial dofs associated with meshes for all benchmark problems.}
\label{table:mesh_hierarchies}
\end{table}

\begin{table}
\centering
\begin{tabular}{r||r|r}
\hline
\multirow{2}{*}{\textbf{Example}}    &  \multicolumn{2}{c}{\textbf{Branch network}}  \\ \cmidrule(lr){2-3} 
                		&  \multicolumn{1}{c|}{\textbf{Layers}} & \textbf{Act.}                               \\ \hline \hline
{Darcy}		& Conv2D[1, 40, 60, 100, 180] + FFN[180, 80, 80, 128] & ReLU	 \\ \hline
{JumpDarcy}		& FFN[1, 256, 256, 256, 128] & Tanh	 \\ \hline
{Wrench/E-shape} 	& FFN[1, 256, 256, 256, 768]  & Tanh	 \\ \hline
{Heat/Snowflake} 	& FFN[1, 256, 256, 256, 128]  & Tanh	 \\ \hline
\hline \hline
\multirow{2}{*}{\textbf{Example}}    &    \multicolumn{2}{c}{\textbf{Trunk network}}    \\ \cmidrule{2-3} 
                		&   \multicolumn{1}{c|}{  \textbf{Layers} } & \textbf{Act.}                                                                        \\ \hline 
{Darcy}		&  FFN[2, 80, 80, 128] &  Tanh \\ \hline
{JumpDarcy}	&  FFN[2, 256, 256, 128] &  Tanh \\ \hline
{Wrench/E-shape} 	&  FFN[3, 512, 256, 256] &  Tanh \\ \hline
{Heat/Snowflake}		&  FFN[3, 256, 256, 128] &  Tanh \\ \hline
\end{tabular}
\caption{The summary of DeepONets' architectures.}
\label{tab:architecture}
\end{table}

\begin{table}[]
\centering
\begin{tabular}{r||r|r}
\hline
\textbf{Example}    		& $\mathbf{N_S}$ 		& \textbf{Time  (mins)}  	\\ \hline \hline
Darcy 				& 2,500				& 198 			   	\\ \hline
JumpDarcy 			& 5,000				& 55 		   			\\ \hline
Wrench 				& 10,000				& 71 					 \\ \hline
E-shape 				& 10,000				& 58 					 \\ \hline
Heat 				& 2,500				& 194				 \\ \hline
Snowflake 			& 2,500				& 77 					\\ \hline
\end{tabular}
\caption{Summary of the number of samples (${N_S}$) and training time for all benchmark problems.}
\label{table:training_params}
\end{table}

\section{Deflation Vectors of NICO Approach}
\label{sec:nico}
For the Darcy problem, we use a constant vector.
For the linear elasticity problem, we use the rigid body modes corresponding to rotation and translation~\cite{jonsthovel2013use}. 
Thus, in three spatial dimensions, for each  $i$-th node of the mesh, the deflation vectors are as 
\begin{equation}
\Pm_i = 
\begin{bmatrix}
1 & 0 & 0 & 0 & z_i & - y_i	 \\
0 & 1 & 0 & -z_i & 0 & x_i	 \\
0 & 0 & 1 & y_i & - x_i & 0	 \\
\end{bmatrix},
\end{equation} 
where $x_i, y_i, z_i$ are coordinates of the $i$-th node.

To construct deflation vectors for the heat equation, we take advantage of the fact that the fully discretized heat equation in two spatial dimensions can be identified with the Helmholtz equation. 
Thus, we consider the following vectors~\cite{vanek1998two}:
\begin{equation}
\Pm_{i} = 
\begin{bmatrix}
e^{-k \frac{\xv_{i} \cdot \dv_{1}}{\Vert \dv_{1} \Vert}} & \cdots & e^{-k \frac{\xv_{i} \cdot \dv_{8}}{\Vert \dv_{8} \Vert}}
\end{bmatrix},
\end{equation}
where $x_{i}, y_{i}$ denote the coordinates of the $i$-th node, and $\{\dv_{j}\}_{j=1}^{8}$ is the set of linearly independent directions given as
${\{\dv_{j}\}_{j=1}^{8}=\{ (1,0), (-1,0), (0,1), (0,-1), (1,1), (-1,-1), (-1,1), (1,-1)\}}$.
 \end{appendices}
 \end{sloppypar}

\bibliographystyle{plain}
\bibliography{biblio}

\end{document}